\documentclass[reqno]{amsart}
\usepackage{amssymb}
\usepackage{pdfsync}
\numberwithin{equation}{section}
\newtheorem{theorem}{Theorem}[section]
\newtheorem{proposition}[theorem]{Proposition}
\newtheorem{corollary}[theorem]{Corollary}

\theoremstyle{definition}

\theoremstyle{definition} %%{remark}t

\newcommand{\bea}{\begin{eqnarray}}
\newcommand{\eea}{\end{eqnarray}}
\newcommand{\beas}{\begin{eqnarray*}}
\newcommand{\eeas}{\end{eqnarray*}}
\newcommand{\beq}{\begin{equation}}
\newcommand{\eeq}{\end{equation}}

\newcommand{\cB}{\mathcal B}

\newcommand{\cD}{\mathcal D}
\newcommand{\cE}{\mathcal E}

\newcommand{\cG}{\mathcal G}

\newcommand{\cR}{\mathcal R}
\newcommand{\cL}{\mathcal L}
\newcommand{\cN}{\mathcal N}

\newcommand{\cMH}{\mathcal{MH}}

\newcommand{\cSM}{\mathcal{SM}}

\newcommand{\BB}{\mathbb B}

\newcommand{\R}{\mathbb R}

\newcommand{\EE}{\mathbb E}
\newcommand{\FF}{\mathbb F}

\newcommand{\vtheta}{\theta}
\newcommand{\Ver}{|\!|}

\newcommand{\ve}{\varepsilon}
\newcommand{\pd}{\partial}

\DeclareMathSymbol{\complement}{\mathord}{AMSa}{"7B}

\def\vv<#1>{\langle #1\rangle}
\def\Vv<#1>{\bigl\langle #1\bigr\rangle}

\begin{document}

% TOPMATTER

\title[Incompressible Two-Phase Flows with Phase Transitions]
{On Well-Posedness of Incompressible\\ Two-Phase Flows with Phase Transitions:\\ The Case of Equal Densities}

\author[J.~Pr\"uss]{J.~Pr\"uss}
\address{Institut f\"ur Mathematik \\
         Martin-Luther-Universit\"at Halle-Witten\-berg\\
         D-60120 Halle, Germany}
\email{jan.pruess@mathematik.uni-halle.de}

\author[Y.~Shibata]{Y.~Shibata}
\address{Department of Mathematics and Research Institute of
         Science and Engineering\\
         JST CREST, Waseda University, Ohkubo 3-4-1, Shinjuku-ku\\
         Tokyo 169-8555, Japan.}
\email{yshibata@waseda.jp}

\author[S.~Shimizu]{S.~Shimizu}
\address{Department of Mathematics, Shizuoka University\\
         Shizuoka 422-8529, Japan}
\email{ssshimi@ipc.shizuoka.ac.jp}

\author[G.~Simonett]{G.~Simonett}
\address{Department of Mathematics\\
         Vanderbilt University \\
         Nashville, TN~37240, USA}
\email{gieri.simonett@vanderbilt.edu}
% \date{}

\thanks{Y.S., S.S., and G.S.\  express their thanks for hospitality
to the Institute of Mathematics, Martin-Luther-Universit\"at Halle-Wittenberg,
where important parts of this work originated.
The research of S.S. was partially supported by challenging Exploratory
Research - 23654048, MEXT, Japan.}

\begin{abstract}
The basic model for incompressible two-phase flows with phase transitions is derived from basic principles and  shown to be thermodynamically consistent in the sense that the total energy is conserved and the total entropy is nondecreasing. The local well-posedness of such problems is proved by means of the technique of
maximal $L_p$-regularity in the case of equal densities. This way we obtain a local semiflow on a well-defined nonlinear state manifold. The equilibria of the system in absence of external forces are identified and it is shown that the negative total entropy is a strict Ljapunov functional for the system. If a solution does not develop singularities, it is proved that it exists globally in time, its orbit is relatively compact, and its limit set is nonempty and contained in the set of equilibria.
\end{abstract}

%%%%%%%%%%%%%%%%%%%%%%
\maketitle

{\small\noindent
\vspace{0.1in}} \vspace{0.1in}\\

\begin{center}
{\small Version of \today}
\end{center}

%\begin{document}

\section{Introduction}
Let $\Omega\subset \R^{n}$ be a bounded domain of class $C^{3-}$, $n\geq2$.
$\Omega$ contains two phases: at time $t$, phase $i$ occupies
the subdomain $\Omega_i(t)$ of
$\Omega$, respectively, with $i=1,2.$
We assume that $\partial \Omega_1(t)\cap\partial \Omega=\emptyset$; this means
that no {\em boundary contact} can occur.
The closed compact hypersurface $\Gamma(t):=\partial \Omega_1(t)\subset \Omega$
forms the interface between the phases.

Let
$u$ denote the velocity field,
$\pi$ the pressure field,
$T$ the stress tensor,
 $\theta>0$  the (absolute) temperature field,
 $\nu_\Gamma$ the outer normal  of $\Omega_1(t)$,
 $V_\Gamma$ the normal velocity of $\Gamma(t)$,
 $H_\Gamma=H(\Gamma(t))=-{\rm div}_\Gamma \nu_\Gamma$ the sum of the principal curvatures of $\Gamma(t)$, and
 $[\![v]\!]=v_2-v_1$ the jump of a quantity $v$ across $\Gamma(t)$.

By  the {\em Incompressible Two-Phase Flow with Phase Transition} we mean the following problem:
 find a family of closed compact hypersurfaces $\{\Gamma(t)\}_{t\geq0}$ contained in $\Omega$
and appropriately smooth functions $u:\R_+\times \bar{\Omega} \to \R^n$, $\pi,\theta:\R_+\times\bar{\Omega}\rightarrow\R$ such that
\begin{equation}
\label{i2pp}
\left\{\begin{aligned}
\partial_t u +u\cdot\nabla u -{\rm div}\, T &=0 && \text{in}&&\Omega\setminus\Gamma(t),\\
T=\mu(\theta)(\nabla u + [\nabla u]^{\sf T}) -\pi I,\quad
{\rm div }\, u &=0 &&\text{in}&&\Omega\setminus\Gamma(t),\\
\kappa (\theta)(\partial_t \theta+u\cdot\nabla\theta)-{\rm div}(d(\theta)\nabla \theta)
-T:\nabla u &=0 &&\text{in}&&\Omega\setminus\Gamma(t),\\
u=\partial_\nu \theta &=0 && \text{on}&&\partial \Omega,\\
[\![u]\!]=[\![\theta]\!]&=0 && \text{on}&& \Gamma(t),\\
[\![T\nu_\Gamma]\!]+\sigma H_\Gamma \nu_\Gamma&=0 &&\text{on}&&\Gamma(t),\\
[\![\psi(\theta)]\!]+\sigma H_\Gamma &=0 &&\text{on}&&\Gamma(t),\\
l(\theta)j+[\![d(\theta)\partial_\nu\theta]\!]&=0 &&\text{on}&&\Gamma(t),\\
V_\Gamma+j- u\cdot \nu_\Gamma&=0 &&\text{on}&&\Gamma(t),\\
 \Gamma(0)=\Gamma_0,\quad u(0,x)=u_0(x), \quad \theta(0,x)&=\theta_0(x) &&\text{in}&&\Omega.
\end{aligned}\right.
\end{equation}
The variable $j$, called {\em phase flux}, can in fact be eliminated from the system: solving for $j$ yields
$$ j=-[\![d(\theta)\partial_\nu\theta]\!]/l(\theta),$$
provided $l(\theta)\neq0$, a property which is assumed later on anyway. In fact, if $l(\theta)=0$ then the interfacial velocity $V$ is not uniquely determined which leads to non-well-posedness of the problem.
The equation for $V_\Gamma$ then becomes
$$V_\Gamma -u\cdot\nu_\Gamma -[\![d(\theta)\partial_\nu\theta]\!]/l(\theta)=0.$$
Note that the sign of the curvature $H_\Gamma$ is negative at a point $x\in\Gamma$ if
$\Omega_1\cap B_r(x)$ is convex, for some sufficiently small $r>0$. Thus if $\Omega_1$ is a ball,
i.e.\ $\Gamma=S_R(x_0)$, then $H_\Gamma=-(n-1)/R$.

Several quantities are derived from the specific {\em Helmholtz free energies} $\psi_i(\theta)$ as follows:
\begin{itemize}
\item
$\epsilon_i(\theta)= \psi_i(\theta)+\theta\eta_i(\theta)$ means the
specific internal energy in phase $i$,
\item
$\eta_i(\theta) =-\psi_i^\prime(\theta)$ the specific entropy,
\item
 $\kappa_i(\theta)= \epsilon^\prime_i(\theta)>0$ the  specific heat capacity, and
\item
$l(\theta)=\theta[\![\psi^\prime(\theta)]\!]=-\theta[\![\eta(\theta)]\!]$ the latent heat.
\end{itemize}
Further $d_i(\theta)>0$ denotes the coefficient of heat conduction in Fourier's law, $\mu_i(\theta)>0$ the viscosity in Newton's law, $\rho:=\rho_1=\rho_2=1$ the constant density,
  and $\sigma>0$ the constant coefficient of surface tension.
 In the sequel we drop the index $i$, as there is no danger of confusion; we just keep in mind that the coefficients depend on the phases.

This model is explained in more detail in the next section. It has been recently proposed by Anderson et al.\ \cite{Gur07}, see also the monographs by Ishii \cite{Ish75} and Ishii and Takashi~\cite{IsTa06}.  We will see below that it is thermodynamically consistent in the sense that in the absence of exterior forces and of fluxes through the outer boundary, the total energy is preserved and the total entropy is nondecreasing. It is in some sense the simplest sharp interface model for incompressible Newtonian two-phase flows taking into account phase transitions driven by temperature.

There is an extensive literature on isothermal incompressible Newtonian two-phase flows without phase transitions, and also on the two-phase Stefan problem with surface tension modeling temperature driven phase transitions. On the other hand, mathematical work on two-phase flow problems including phase transitions is rare. In this direction, we only know of the papers by Hoffmann and Starovoitov \cite{HoSt98a,HoSt98b}, dealing with a simplified two-phase flow model, and Kusaka and Tani \cite{KuTa99,KuTa02} which are two-phase for the
temperature but with only one phase moving. The papers of DiBenedetto and Friedman \cite{DBFr86} and DiBenedetto and O'Leary\cite{DBOL93} deal with weak solutions of conduction-convection problems with phase change. However,
the models considered in these papers do not seem to be consistent with thermodynamics.

It is the purpose of this paper to present a rigorous analysis of problem \eqref{i2pp} in the framework of $L_p$-theory. We prove local well-posedness, and we identify the equilibria of the problem.
In turns out that the equilibria are the same as those for the thermodynamically consistent two-phase Stefan problem with surface tension. This heavily depends on the fact that the densities of the two phases are assumed to be equal; in this case the problem is {\em temperature dominated}.

In a forthcoming paper we will consider the case where the densities are not equal; then the
solution behavior is different, as the interfacial mass flux has a direct impact on the velocity field of the fluid. Physically this relates to the so-called {\em Stefan currents} which are induced by phase transitions.
The velocity field is then no longer continuous across the interface and this leads to different analytic properties of the model. We call this case {\em velocity dominated}.

The plan for this paper is as follows. In the next section we introduce and discuss the model in some detail and show that it is thermodynamically consistent. In Section 3 we prove that the negative total entropy is a strict Ljapunov functional for the problem, i.e.\ is strictly decreasing along nonconstant solutions. We identify the equilibria which consist of velocity zero, constant temperature, constant pressures in the phases, and the disperse phase $\Omega_1$ consists of finitely many nonintersecting balls of the {\em same radius}. This is in contrast to the isothermal two-phase Navier-Stokes equation with surface tension where phase transitions are neglected: in that case the balls may have different sizes!

In our approach we apply the well-known {\em direct mapping method} where the problem with moving interface is transformed to a problem on a fixed domain. This transformation will be introduced in Section~4; as a result we  obtain a quasilinear parabolic evolution problem with dynamic boundary condition. In Section~5 we study  maximal $L_p$-regularity of the underlying linearized problem. It turns out that the linear problem is equivalent to the linearized Stefan problem with surface tension, followed by a Stokes problem where the curvature serves as an input variable on the interface. Here we rely on previous work Pr\"uss, Simonett and Zacher \cite{PSZ10}. Local well-posedness is discussed in Section~6, while the proof of the main result is given in Section~7. It is based on maximal $L_p$-regularity of the linear problem and the contraction mapping principle.
Finally, in Section~8 we consider global existence of solutions. It is shown that
solutions exist globally in time, provided the interface satisfies a uniform ball condition,
the temperature and the modulus of the latent heat are bounded from below, and certain a priori estimates on the solution hold. In a forthcoming paper we will show that in this situation the solution converges
to an equilibrium in the topology of the state manifold.

\section{The Model}
\noindent
In this section we explain the model.
We begin with

\medskip

\noindent
{\bf Balance of Mass:}
\begin{equation*}
\begin{aligned}
\partial_t \rho + {\rm div}\, (\rho u)&=0 &&\text{in}&& \Omega\setminus \Gamma(t),\\
[\![\rho(u-u_\Gamma)]\!]\cdot\nu_\Gamma &=0 &&\text{on}&& \Gamma(t),
%\noindent
\end{aligned}
\end{equation*}
where  $u_\Gamma$ is the interfacial velocity.
Hence the normal velocity $V_\Gamma$ can be expressed by $V_\Gamma=u_\Gamma\cdot\nu_\Gamma$.
We define the {\em phase flux}, more precisely, the {\em interfacial mass flux}, by means of
\begin{equation}
 j:= \rho(u-u_\Gamma)\cdot\nu_\Gamma,
 \quad \mbox{i.e.}\quad [\![\frac{1}{\rho}]\!] j= [\![u\cdot \nu_\Gamma]\!],
\end{equation}
and we note that $j$ is independent of the phases, hence well-defined.
Therefore, a {\em phase transition} takes place if $j\neq 0$.
If $j\equiv0$ then $u\cdot\nu_\Gamma=u_\Gamma\cdot\nu_\Gamma=V_\Gamma$
and in this case the interface is advected with the velocity field $u$.

In this paper we consider the {\em completely incompressible} case,
i.e.\ we assume that the densities are constant in the phases $\Omega_i$.
Then conservation of mass reduces to
\begin{equation}
\label{incompressible}
 {\rm div}\, u=0\quad\text{in}\quad \Omega\setminus \Gamma(t).
\end{equation}
If only  the latter property holds, we say that the material is {\em incompressible}.

Integrating over the finite domain $\Omega_1(t)$ we obtain by the Reynolds transport theorem
\begin{equation*}
\begin{aligned}
\frac{d}{dt} \int_{\Omega_1(t)}\rho\, dx &= \int_{\Gamma(t)} \rho V_\Gamma\,ds
+\int_{\Omega_1(t)} \partial_t\rho\, dx \\
&=\int_{\Gamma(t)} \rho V_\Gamma\, ds -\int_{\Omega_1(t)} {\rm div}(\rho u)\, dx \\
&= \int_{\Gamma(t)} (\rho u_\Gamma\cdot\nu_\Gamma -\rho u\cdot \nu_\Gamma)\,ds
= -\int_{\Gamma(t)} j\, ds.
\end{aligned}
\end{equation*}
In case $|\Omega_2(t)|$ is finite as well, we obtain in the same way
$$\frac{d}{dt} \int_{\Omega_2(t)}\rho\, dx = \int_{\Gamma(t)} j\, ds,$$
proving conservation of total mass.

Therefore, if both phases are completely incompressible  we have
$$\rho_1|\Omega_1(t)|+\rho_2|\Omega_2(t)|\equiv \rho_1|\Omega_1(0)|+\rho_2|\Omega_2(0)|=:c_0.$$
This implies
$$[\![\rho]\!]|\Omega_1(t)|=\rho_2|\Omega|-c_0,$$
hence $|\Omega_1(t)|=constant$ in the case of nonequal densities, i.e.\ the phase volumes are preserved.
If the densities are equal, there is no restriction on the phase volumes due to conservation of mass.
This constitutes a big difference betweeen the two cases.

In later sections of this paper we are concerned with the completely incompressible case
where, in addition, the densities are equal, and hence conservation of mass does then not imply conservation of phase volumes.

Next we consider
\medskip

\noindent
{\bf Balance of Momentum:}
\begin{equation*}
\begin{aligned}
\partial_t (\rho u) + {\rm div}\, (\rho u\otimes u) -{\rm div}\,T&=\rho f
&&\text{in}&& \Omega\setminus \Gamma(t),\\
[\![\rho u\otimes(u-u_\Gamma)-T]\!]\nu_\Gamma &={\rm div}_\Gamma T_\Gamma &&\text{on}&&\Gamma(t).
\end{aligned}
\end{equation*}
Using balance of mass and the definition of the phase flux $j$ we may rewrite this conservation
law as follows:
\begin{equation}
\label{momentum}
\begin{aligned}
\rho(\partial_t u + u\cdot\nabla u) -{\rm div}\,T&=\rho f &&\text{in}&& \Omega\setminus \Gamma(t),\\
[\![u]\!]j-[\![T\nu_\Gamma]\!] &={\rm div}_\Gamma T_\Gamma &&\text{on}&& \Gamma(t).
\end{aligned}
\end{equation}
Here $T_\Gamma$ denotes surface stress, which in this paper will be assumed to be of the form
\begin{equation*}
 T_\Gamma = \sigma P_\Gamma,\quad P_\Gamma=I -\nu_\Gamma\otimes\nu_\Gamma,
\end{equation*}
where $\sigma>0$ is a constant. So in particular we neglect here surface viscosities as well as Marangoni forces.
Then we have
\begin{equation}
\label{surface-stress}
{\rm div}_\Gamma T_\Gamma= \sigma H_\Gamma \nu_\Gamma
\end{equation}
which is the usual assumption for surface stress.
\smallskip

We now turn our attention to

\medskip

\noindent
{\bf Balance of Energy:}\\
Let $q$ denote the heat flux and $r$ an external heat source. Then balance of energy becomes
\begin{equation*}
\begin{aligned}
\partial_t (\frac{\rho}{2} |u|_2^2 +\rho \epsilon) + {\rm div}\,\{ (\frac{\rho}{2}  |u|_2^2 + \rho \epsilon) u\}
-{\rm div} (Tu-q)&=\rho f\cdot u +\rho r &&\text{in}&&\Omega\setminus \Gamma(t),\\
[\![(\frac{\rho}{2} |u|_2^2+\rho \epsilon)(u-u_\Gamma)-(Tu-q)]\!]\cdot\nu_\Gamma
&={\rm div}_\Gamma T_\Gamma\cdot u_\Gamma && \text{on} &&  \Gamma(t).
\end{aligned}
\end{equation*}
Note that the second line implies that the only surface energy taken into account here is that induced by surface tension.

Using balance of mass, balance of momentum, and the definition of the phase flux $j$ we may rewrite
this conservation law as follows:
\begin{equation}
\begin{aligned}
\rho(\partial_t \epsilon + u\cdot\nabla \epsilon) +{\rm div}\,q-T:\nabla u&=\rho r
&&\text{in}&& \Omega\setminus \Gamma(t),\\
([\![\epsilon]\!]+[\![\frac{1}{2}|u-u_\Gamma|_2^2]\!])j-[\![T\nu_\Gamma\cdot(u-u_\Gamma)]\!] + [\![q\cdot\nu_\Gamma]\!] &=0 &&\text{on}&&\Gamma(t).
\end{aligned}
\end{equation}
\goodbreak
The total energy is given by
\begin{equation*}
\label{energy}
{\sf E}:={\sf E}(u,\theta,\Gamma):=\frac{1}{2}\int_\Omega \rho|u|_2^2\,dx
+ \int_\Omega \rho \epsilon(\theta)\,dx +\sigma |\Gamma|,
\end{equation*}
where $|\Gamma|$ denotes the surface area of $\Gamma$.
Note that here the assumption of $\sigma$ being constant is important! Otherwise, one has to take into account surface energy, {\em locally}, which means that we would have to include a balance of surface energy.

In the absence of external forces $f$ and heat sources $r$, for the time derivative of ${\sf E}$ we obtain by the transport theorem
\begin{align*}
\frac{d}{dt}{\sf E} &= \int_\Omega \{u\cdot\rho\partial_t u + \rho\partial_t \epsilon(\theta)\}\,dx
-\int_\Gamma \{[\![\frac{\rho}{2}|u|_2^2 +\rho \epsilon(\theta)]\!]+\sigma H_\Gamma\}V_\Gamma\,ds\\
&= -\int_\Omega  \{\rho(u\cdot\nabla)u\cdot u -{\rm div}\, T \cdot u
+\rho (u\cdot\nabla)\epsilon(\theta) +{\rm div}\, q - T:\nabla u\}\,dx\\
&\quad-\int_\Gamma\big \{[\![\frac{\rho}{2}|u|_2^2 +\rho \epsilon(\theta)]\!]
+\sigma H_\Gamma\big\}u_\Gamma\cdot\nu_\Gamma\, ds\\
&= \int_\Gamma\big\{[\![(\frac{\rho}{2}|u|_2^2 +\rho \epsilon(\theta))(u-u_\Gamma)\cdot\nu_\Gamma ]\!]
-[\![T\nu_\Gamma\cdot u ]\! ]+[\![q\cdot\nu_\Gamma]\!]-\sigma H_\Gamma V_\Gamma\big\}\,ds\\
&=\int_\Gamma\big\{[\![\frac{1}{2}|u|_2^2 + \epsilon(\theta) ]\!]j-[\![Tu\cdot \nu_\Gamma ]\! ]+[\![q\cdot\nu_\Gamma]\!]-\sigma H_\Gamma V_\Gamma\big\}\,ds,
\end{align*}
provided that energy transport through the outer boundary is zero, which means
$$\big\{\frac{\rho|u|^2}{2}+\rho\epsilon(\theta)\big\}u\cdot\nu -T\nu\cdot u +q\cdot\nu=0\quad \text{on} \quad \partial\Omega.$$
Hence
$$ \frac{d}{dt} {\sf E}(u,\theta,\Gamma)=0,$$
which implies that the total energy is preserved.

\medskip

\noindent
{\bf Entropy Production:}\\
We introduce now the fundamental thermodynamic relations which read
\begin{equation}
\begin{aligned}
\label{thermodynamics}
\epsilon(\theta)&=\psi(\theta)+\theta\eta(\theta),\quad \eta(\theta)=-\psi^\prime(\theta),\\
\kappa(\theta)&=\epsilon^\prime(\theta),
\hspace{1.7cm} l(\theta)=\theta[\![\psi^\prime(\theta)]\!]=-\theta[\![\eta(\theta)]\!],\nonumber
\end{aligned}
\end{equation}
where $\psi(\theta)$ means the {\em Helmholtz free energy} which should be considered as given, but depends on the phases. The quantities $\eta$ and $\kappa$, $l$ are called {\em entropy}, {\em heat capacity}, and  {\em latent heat}, respectively.

The total entropy is defined by
\begin{equation*}
\label{entropy}
\Phi(\theta,\Gamma)=\int_\Omega \rho\eta(\theta)\,dx.
\end{equation*}
The main idea of the modeling approach presented here is that there should be no entropy production on the interface, as it is considered to be ideal; in particular it is assumed to carry no mass and no energy except for surface tension.
We have with \eqref{incompressible}-\eqref{constitutive-interface} for $f=r=0$,
and with $\epsilon^\prime(\theta)=\theta\eta^\prime(\theta)$,
\begin{align*}
\frac{d}{dt}& \Phi(u,\theta,\Gamma) = \int_\Omega \rho \partial_t \eta(\theta)\,dx
 - \int_\Gamma [\![\rho \eta(\theta)]\!] V_\Gamma\, ds\\
&= \int_\Omega\rho\eta^\prime(\theta)\partial_t\theta\,dx
- \int_\Gamma [\![\rho \eta(\theta)]\!] u_\Gamma\cdot\nu_\Gamma\,ds\\
&= \int_\Omega \Big\{\frac{\eta^\prime(\theta)}{\epsilon^\prime(\theta)}
\{ T:D -{\rm div}\, q\}-\rho u\cdot\nabla \eta(\theta)\Big\}\,dx
- \int_\Gamma [\![\rho \eta(\theta)]\!] u_\Gamma\cdot\nu_\Gamma\, ds\\
&= \int_\Omega\Big\{\frac{T:D}{\theta} -\frac{q\cdot\nabla\theta}{\theta^2}\Big\}\,dx
+\int_\Gamma\Big\{[\![\rho\eta(\theta)(u-u_\Gamma)\cdot\nu_\Gamma]\!]
+\frac{1}\theta [\![q\cdot\nu_\Gamma]\!]\Big\}\,ds \\
&=\int_\Omega\Big\{\frac{T:D}{\theta} -\frac{q\cdot\nabla\theta}{\theta^2}\Big\}\,dx
+ \int_\Gamma\frac{1}{\theta}\Big\{[\![-l(\theta)j + [\![q\cdot\nu_\Gamma]\!]\Big\}\, ds,
\end{align*}
provided there is no entropy flux through the outer boundary $\partial\Omega$, which means
$$q\cdot\nu + \rho\theta\eta(\theta) u\cdot\nu=0 \quad \text{on}\quad \partial\Omega.$$
To ensure conservation of energy as well as entropy through the outer boundary we impose in this paper for simplicity the following

\medskip

\noindent
{\bf Constitutive Laws on the Outer Boundary $\partial\Omega$:  }
\begin{equation}
\label{constitutive-outer}
q\cdot\nu=0,\quad u=0.
\end{equation}
Actually, we may consider more general conditions at the outer boundary, like $q\cdot\nu=0$ and the partial slip condition $u\cdot\nu=0$, $P_{\partial\Omega} T\nu=0$, but we refrain from doing this here.

As constitutive laws in the phases we employ Newton's law for the stress tensor and Fourier's law for the heat flux; these  ensure nonnegative entropy production in the bulk. Recall that  we consider the completely incompressible case.

\medskip

\noindent
{\bf Constitutive Laws in the Phases:}
\begin{equation}
\label{constitutive-phases}
\begin{aligned}
T&=S-\pi I,\quad
S=2\mu(\theta) D,\quad D = \frac{1}{2}(\nabla u +[\nabla u]^{\sf T}),\\
q& = -d(\theta) \nabla\theta.
\end{aligned}
\end{equation}
We assume $\kappa(\theta)=\epsilon^\prime(\theta)=-\theta\psi^{\prime\prime}(\theta)>0$, as well  $\mu(\theta)>0$ for the shear viscosity, and $d(\theta)>0$  for the heat diffusivity; note that each of these quantities may jump across the interface.

Finally, on the interface we assume
\medskip
\goodbreak
\noindent
{\bf Constitutive Laws on the Interface $\Gamma(t)$}:
\begin{equation}\label{constitutive-interface}
 [\![\theta]\!]=[\![P_\Gamma u]\!] =0,\quad -l(\theta)j+[\![q\cdot\nu_\Gamma]\!]=0.
\end{equation}
Hence in this model the temperature and the tangential part of the velocity are continuous across the interface; the latter means that there is no tangential slip at the interface, as it is considered ideal. The third equation in \eqref{constitutive-interface} means that no entropy is generated on the interface. This implies the relation
\begin{equation}\label{entropy-production}
\frac{d}{dt}\Phi(\theta,\Gamma)=
\int_\Omega \Big\{\frac{2\mu(\theta)}{\theta}|D|_2^2 +\frac{d(\theta)}{\theta^2}|\nabla\theta|_2^2\Big\}\,dx\geq0
\end{equation}
for the entropy production.\\
Employing the constitutive assumption $[\![P_\Gamma u]\!]=0$ with
$P_\Gamma[\![T\nu_\Gamma]\!]=0$ yields
\begin{equation*}
\begin{split}
&[\![|u-u_\Gamma|^2]\!]= [\![\frac{1}{\rho^2}]\!]j^2,\\
&[\![T\nu_\Gamma\cdot(u-u_\Gamma)]\!]
=[\![T\nu_\Gamma\cdot P_\Gamma(u-u_\Gamma)]\!]+[\![\frac{1}{\rho}T\nu_\Gamma\cdot\nu_\Gamma]\!]j
= [\![\frac{1}{\rho}T\nu_\Gamma\cdot\nu_\Gamma ]\!]j.
\end{split}
\end{equation*}
Combining these conditions with the energy balance across the interface and the constitutive law for $q$
results into the following jump conditions, assuming that $j$ may take arbitrary values:
\begin{equation}
\begin{aligned}
l(\theta)j +[\![d(\theta)\partial_\nu\theta]\!]&=0 &&\text{on}&&\Gamma(t),\\
\{[\![\psi(\theta)]\!]+[\![\frac{1}{2\rho^2}]\!]j^2 -[\![\frac{1}{\rho}T\nu_\Gamma\cdot\nu_\Gamma]\!]\}j&=0
&&\text{on}&&\Gamma(t).
\end{aligned}
\end{equation}
The first of these two equations is the {\em Stefan law} and the second the (generalized) {\em Gibbs-Thomson law}.

Note that if $j\equiv 0$, i.e.\ in the absence of a phase transition,
%and the densities are equal,
the Stefan law becomes $[\![d(\theta)\partial_\nu\theta]\!]=0$,
and the Gibbs-Thomson relation trivializes. On the other hand, if phase transitions occur then $j$ must be considered as arbitrary, and so the Gibbs-Thomson law becomes
\begin{equation}
[\![\psi(\theta)]\!]+[\![\frac{1}{2\rho^2}]\!]j^2 -[\![\frac{1}{\rho}T\nu_\Gamma\cdot\nu_\Gamma]\!]=0.
\end{equation}
We summarize to obtain the following model for incompressible two-phase flows with phase transitions.
\goodbreak
\bigskip

\noindent
{\bf The Complete Model}
\vspace{2mm}\\
{\em In the bulk} $\Omega\setminus\Gamma(t)$:
\begin{equation}
\label{bulk}
\begin{aligned}
\rho(\partial_t u + u\cdot\nabla u) -{\rm div}\,T=0,&&\\
T=\mu(\theta)(\nabla u+[\nabla u]^{\sf T})-\pi I,\quad
{\rm div}\, u=0,&&\\
\rho\kappa(\theta)(\partial_t \theta  + u\cdot\nabla \theta) -{\rm div}(d(\theta)\nabla\theta)-T:\nabla u =0.&&
\end{aligned}
\end{equation}
\vspace{2mm}\\
{\em On the interface} $\Gamma(t)$:
\begin{equation}
\label{interface}
\begin{aligned}
\mbox{}[\![\frac{1}{\rho}]\!]j^2\nu_\Gamma-[\![T\nu_\Gamma]\!]-\sigma H_\Gamma \nu_\Gamma=0,&&
\quad\mbox{}[\![u]\!]&= [\![\frac{1}{\rho}]\!] j \nu_\Gamma,\\
\quad l(\theta)j + [\![d(\theta)\partial_\nu \theta]\!]=0,&&
 [\![\theta]\!]&=0,\\
\mbox{}[\![\psi(\theta)]\!]+[\![\frac{1}{2\rho^2}]\!]j^2 -[\![\frac{1}{\rho}T\nu_\Gamma\cdot\nu_\Gamma]\!]=0,&&
\quad V_\Gamma &= u\cdot\nu_\Gamma -\frac{1}{\rho} j.
\end{aligned}
\end{equation}
\vspace{2mm}\\
{\em On the outer boundary} $\partial\Omega$:
\begin{equation}
\label{BC}
u=0,\quad\partial_\nu \theta =0.
\end{equation}
{\em Initial conditions}:
\begin{equation}
\label{IC}
\Gamma(0)=\Gamma_0,\quad u(0)=u_0,\quad \theta(0)=\theta_0.
\end{equation}
This results in the model \eqref{i2pp} if we assume in addition that the densities are equal, i.e.\ $\rho_1=\rho_2=1$.

\section{Equilibria}
As we have seen in Section 2, the negative total entropy is a Ljapunov functional for the problem,
and it is even a strict one. To see this, assume that $\Phi$ is constant on some interval $(t_1,t_2)$.
Then $\frac{d}{dt}\Phi(\theta,\Gamma)=0$ in $(t_1,t_2)$, hence
$D=0$ and $\nabla\theta=0$ in $(t_1,t_2)\times\Omega$.
Therefore,
$\theta$ is constant which implies $[\![d(\theta)\partial_\nu\theta]\!]=0$, and then from the interfacial boundary condition we obtain $j=0$, provided $l(\theta)\neq0$.
 Then $[\![u]\!]=0$, hence by Korn's inequality we have $\nabla u=0$ and then $u=0$ by the no-slip condition
 on $\partial\Omega$.
 This implies further $\partial_t\theta=\partial_tu =0$ and $V_\Gamma=0$, i.e.\ we are at equilibrium.
 Further, $\nabla\pi=0$, i.e.\ the pressure is constant in the components of the phases, and
 $\sigma H_\Gamma=[\![\pi]\!]$, $[\![\psi(\theta)]\!]= -[\![\pi]\!]$ are constant as well. Since $\theta$ is continuous across the interface the last relation shows that $\pi$ is constant in all of $\Omega_1$, even if it is not connected.
 From this we finally deduce that
  $\Omega_1$ is a ball if it is connected, or a finite union of balls of equal radii, as $\Omega$ is bounded.
Let us summarize.
\goodbreak
%%%%%%%%%%%%%%%%%%%%%%
\begin{theorem}
Let $\sigma>0$, $\psi_i\in C^3(0,\infty)$, $\mu_i,d_i\in C^2(0,\infty)$,
 $$-\psi^{\prime\prime}_i(s),\mu_i(s),d_i(s)>0, \quad [\![\psi^\prime(s)]\!]\neq0,$$
 for $s>0$, $i=1,2$. Then the following assertions hold.
\begin{itemize}
\item
[\bf (i)] The total energy ${\sf E}$ is conserved along smooth solutions.
\vspace{2mm}
\item
[\bf (ii)] The negative total entropy $-\Phi$ is a strict Ljapunov functional, which means that $-\Phi$ is strictly decreasing along nonconstant smooth solutions.
\vspace{2mm}
\item
[\bf (iii)] The non-degenerate equilibria (i.e.\ no boundary contact) are zero velocity, constant temperature, constant pressure in each phase, and $\Omega_1$ consists of a finite number of nonintersecting balls of equal size.
\end{itemize}
\end{theorem}

Therefore the equilibria are the same as those for the thermodynamically consistent Stefan problem with surface tension which  has been discussed recently in Pr\"uss, Simonett and Zacher \cite{PSZ10}.
We now discuss some basic considerations from that paper.
Suppose $(u,\pi,\theta,\Gamma)$ is an equilibrium with
\begin{equation*}
 \Gamma =\Big\{\bigcup_{j=1}^m S(x_j,R): B(x_j,R)\cap B(x_k,R)=\emptyset,\; j\neq k,
 \ B(x_j,R)\subset\Omega\Big\},
\end{equation*}
where $B(x,R)$ denotes the ball with center $x$ and radius $R$, and $S(x,R)$ its boundary.
To determine $([\![\pi]\!],\theta, R)$ we have to solve the system
\begin{align}
\label{eq-rel}
\varphi(\theta,R):= |\Omega_1| \epsilon_1(\theta) +|\Omega_2| \epsilon_2(\theta)
+ \sigma|\Gamma|&={\sl E}_0,\nonumber\\
[\![\pi]\!]&=\sigma H_\Gamma,\\
[\![\psi(\theta)]\!]&=-[\![\pi]\!],\nonumber
\end{align}
where  ${\sf E}_0={\sf E}(u_0,\theta_0,\Gamma_0)$ is the prescribed initial total energy, and $\varphi(\theta,R)$ means the total energy at a given equilibrium to be studied.
We have
$$[\![\pi]\!]=\sigma H_\Gamma=-(n-1)\sigma/R,\quad R=R(\theta)=(n-1)\sigma/[\![\psi(\theta)]\!].$$
Hence there remains a single equation for the
equilibrium temperature $\theta$, namely
$${\sf E}_0 =  \varphi(\theta):=|\Omega|\epsilon_2(\theta) -m(\omega_n/n) R^n(\theta)[\![\epsilon(\theta)]\!]+ \sigma m \omega_n R^{n-1}(\theta),$$
where $\omega_n$ denotes the area of the unit sphere in $\R^n$.
Note that only the temperature range $[\![\psi(\theta)]\!]>0$ is relevant due to the requirement $R>0$, and with
$$R_m^*=\sup\{R>0:\, \Omega \mbox{ contains } m \mbox{ disjoint balls of radius } R\}$$
we must also have $R<R_m^*$, i.e. with  $h(\theta) = [\![\psi(\theta)]\!]$
$$h(\theta)>  \frac{\sigma(n-1)}{R_m^*}.$$
With $\epsilon(\theta)=\psi(\theta)-\theta \psi^\prime(\theta)$ we may rewrite $\varphi(\theta)$ as
$$ \varphi(\theta) =|\Omega|\epsilon_2(\theta) + c_n \Big( \frac{1}{h(\theta)^{n-1}} +(n-1)\theta\frac{h^\prime(\theta)}{h(\theta)^n}\Big),$$
where we have set $\displaystyle c_n = m\frac{\omega_n}{n(n-1)}((n-1)\sigma)^n.$

Next, with
$$ R^\prime(\theta) = - \frac{\sigma(n-1)h^\prime(\theta)}{ h^2(\theta)} = - \frac{ h^\prime(\theta) R^2(\theta)}{\sigma(n-1)}$$
we obtain
\begin{align*}
\varphi^\prime(\theta)&= |\Omega|\epsilon^\prime_2(\theta) -[\![\epsilon^\prime(\theta)]\!]|\Omega_1|
+ m\omega_n\Big(\frac{\sigma(n-1)}{ R(\theta)} -[\![e(\theta)]\!]\Big)R^{n-1}(\theta)R^{\prime}(\theta)\\
&=|\Omega|\kappa_2(\theta) -[\![\kappa(\theta)]\!]|\Omega_1|+ m\omega_n \theta h^\prime(\theta) R^{n-1}(\theta) R^\prime(\theta)\\
&= (\kappa|1)_\Omega -\theta h^\prime(\theta)|\Gamma| \frac{ h^\prime(\theta) R^2(\theta)}{\sigma(n-1)}\\
&= \frac{(\kappa|1)_\Omega R^2(\theta)}{\sigma(n-1)}\Big\{ \frac{\sigma(n-1)}{ R^2(\theta)} - \frac{l(\theta)^2|\Gamma|}{\theta(\kappa|1)_\Omega}\Big\},
\end{align*}
with $l(\theta)=\theta h^\prime(\theta)$. It will turn out that the term in the parentheses determines whether an equilibrium is stable:
it is stable if and only if  $m=1$ and $\varphi^\prime(\theta)<0$. In particular, if $\Omega_1$ is not connected this equilibrium
is unstable; this fact is related to the physical phenomenon called {\em Ostwald ripening}.

In general it is not a simple task to analyze the equation for the temperature
$$\varphi(\theta)= |\Omega|\epsilon_2(\theta) + c_n \Big( \frac{1}{h(\theta)^{n-1}} +(n-1)\theta\frac{h^\prime(\theta)}{h(\theta)^n}\Big)={\sf E}_0,$$
unless more properties of the functions $\epsilon_2(\theta)$ and, in particular, of $[\![\psi(\theta)]\!]$ are known; cf.\ Pr\"uss, Simonett and Zacher \cite{PSZ10} for further discussion and results.

\section{Transformation to a Fixed Domain}

Let $\Omega\subset\R^n$ be a bounded domain with boundary $\partial\Omega$ of class $C^2$, and suppose
$\Gamma\subset\Omega$ is a closed hypersurface of class $C^2$,
i.e.\ a $C^2$-manifold which is the boundary of a bounded domain
$\Omega_1\subset\Omega$; we then set
$\Omega_2=\Omega\backslash\bar{\Omega}_1$. Note that $\Omega_2$ is connected,
but $\Omega_1$ maybe  disconnected, however, it consists of finitely
many components only, since $\partial\Omega_1=\Gamma$ by assumption is a manifold, at least  of class $C^2$.
Recall that the {\em second order bundle} of $\Gamma$ is given by
$$\cN^2\Gamma:=\{(p,\nu_\Gamma(p),\nabla_\Gamma\nu_\Gamma(p)):\, p\in\Gamma\}.$$
Here $\nu_\Gamma(p)$ denotes the outer normal of $\Omega_1$ at $p\in\Gamma$ and $\nabla_\Gamma$ the surface gradient on $\Gamma$. The Weingarten tensor $L_\Gamma$ on $\Gamma$ is defined by
$$L_\Gamma(p)= -\nabla_\Gamma\nu_\Gamma(p),\quad p\in\Gamma.$$
The eigenvalues $\kappa_j(p)$ of $L_\Gamma(p)$ are the principal curvatures of $\Gamma$ at $p\in\Gamma$, and we have
$|L_\Gamma(p)|=\max_j|\kappa_j(p)|$. The curvature $H_\Gamma(p)$ (more precisely $(n-1)$ times mean curvature) is defined as the trace of $L_\Gamma(p)$, i.e.
$$ H_\Gamma(p)= \sum_{j=1}^{n-1} \kappa_j(p) = -{\rm div}_\Gamma \nu_\Gamma(p),$$
where ${\rm div}_\Gamma$ means surface divergence.
Recall also the
{\em Hausdorff distance} $d_H$ between the two closed subsets $A,B\subset\R^m$,
defined by
$$d_H(A,B):= \max\{\sup_{a\in A}{\rm dist}(a,B),\sup_{b\in B}{\rm dist}(b,A)\}.$$
Then we may approximate $\Gamma$ by a real analytic hypersurface $\Sigma$ (or merely $\Sigma\in C^3$),
in the sense that the Hausdorff distance of the second order bundles of
$\Gamma$ and $\Sigma$ is as small as we want. More precisely, for each $\eta>0$ there is a
real analytic closed hypersurface $\Sigma$ such that
$d_H(\cN^2\Sigma,\cN^2\Gamma)\leq\eta$. If $\eta>0$ is small enough, then $\Sigma$
bounds a domain $\Omega_1^\Sigma$ with  $\overline{\Omega_1^\Sigma}\subset\Omega$, and we set $\Omega_2^\Sigma=\Omega\backslash\overline{\Omega_1^\Sigma}$.

It is well known that a hypersurface $\Sigma$ of class $C^2$ admits a tubular neighborhood,
which means that there is $a>0$ such that the map
\begin{eqnarray*}
&&\Lambda: \Sigma \times (-a,a)\to \R^n \\
&&\Lambda(p,r):= p+r\nu_\Sigma(p)
\end{eqnarray*}
is a diffeomorphism from $\Sigma \times (-a,a)$
onto $\cR(\Lambda)$. The inverse
$$\Lambda^{-1}:\cR(\Lambda)\mapsto \Sigma\times (-a,a)$$ of this map
is conveniently decomposed as
$$\Lambda^{-1}(x)=(\Pi_\Sigma(x),d_\Sigma(x)),\quad x\in\cR(\Lambda).$$
Here $\Pi_\Sigma(x)$ means the orthogonal projection of $x$ to $\Sigma$ and $d_\Sigma(x)$ the signed
distance from $x$ to $\Sigma$; so $|d_\Sigma(x)|={\rm dist}(x,\Sigma)$ and $d_\Sigma(x)<0$ if and only if
$x\in \Omega_1^\Sigma$. In particular we have $\cR(\Lambda)=\{x\in \R^n:\, {\rm dist}(x,\Sigma)<a\}$.

Note that on the one hand, $a$ is determined by the curvatures of $\Sigma$, i.e.\ we must have
$$0<a<\min\{1/|\kappa_j(p)|: j=1,\ldots,n-1,\; p\in\Sigma\},$$
where $\kappa_j(p)$ mean the principal curvatures of $\Sigma$ at $p\in\Sigma$.
But on the other hand, $a$ is also connected to the topology of $\Sigma$,
which can be expressed as follows. Since $\Sigma$ is a compact $C^2$ manifold of
dimension $n-1$ it satisfies the ball condition, which means that
there is a radius $r_\Sigma>0$ such that for each point $p\in \Sigma$
there are points $x_j\in \Omega_i^\Sigma$, $i=1,2$, such that $B_{r_\Sigma}(x_j)\subset \Omega_i^\Sigma$, and
$\bar{B}_{r_\Sigma}(x_i)\cap\Sigma=\{p\}$. Choosing $r_\Sigma$ maximal,
we then must also have $a<r_\Sigma$. In the sequel we fix
$$a=\frac{1}{2} \min\{ r_\Sigma,\frac{1}{|\kappa_j(p)|}: j=1,\ldots, n-1, p\in\Sigma\}.$$
For later use we note that the derivatives of $\Pi_\Sigma(x)$ and $d_\Sigma(x)$ are given by
$$\nabla d_\Sigma(x)= \nu_\Sigma(\Pi_\Sigma(x)),
\quad D\Pi_\Sigma(x) = M_0(d_\Sigma(x))(\Pi(x))P_\Sigma(\Pi_\Sigma(x))$$
for $|d_\Sigma(x)|<a$,
where $P_\Sigma(p)=I-\nu_\Sigma(p)\otimes\nu_\Sigma(p)$ denotes the orthogonal projection onto
the tangent space $T_p\Sigma$ of $\Sigma$ at $p\in\Sigma$, and
\begin{equation*}
M_0(r)(p)=(I-r L_\Sigma(p))^{-1},\quad (r,p)\in (-a,a)\times\Sigma.
\end{equation*}
Note that $$|M_0(r)(p)|\leq 1/(1-r|L_\Sigma(p)|)
\leq 2,\quad \mbox{ for all } (r,p)\in (-a,a)\times\Sigma.$$
Setting $\Gamma=\Gamma(t)$, we may use the map $\Lambda$ to parameterize the unknown free
boundary $\Gamma(t)$ over $\Sigma$ by means of a height function $h(t,p)$ via
$$\Gamma(t)=\{ p+ h(t,p)\nu_\Sigma(p): p\in\Sigma,\; t\geq0\},$$
at least for small $|h|_\infty$.
Extend this diffeomorphism to all of $\bar{\Omega}$ by means of
$$ \Xi_h(t,x) = x +\chi(d_\Sigma(x)/a)h(t,\Pi_\Sigma(x))\nu_\Sigma(\Pi_\Sigma(x))=:x+\xi_h(t,x).$$
Here $\chi$ denotes a suitable cut-off function. More precisely, $\chi\in\cD(\R)$,
$0\leq\chi\leq 1$, $\chi(r)=1$ for $|r|<1/3$, and $\chi(r)=0$ for $|r|>2/3$.
 Note that $\Xi_h(t,x)=x$ for $|d(x)|>2a/3$, and
%$$\Xi_h^{-1}(t,x)= x-h(t,\Pi_\Sigma(x))\nu_\Sigma(\Pi_\Sigma(x))\quad \mbox{ for }\; |d_\Sigma(x)|<a/3,$$
%in particular,
$$\Xi_h^{-1}(t,x)= x-h(t,x)\nu_\Sigma(x),\quad x\in\Sigma,$$
for $|h|_\infty$ sufficiently small.
Now we define the transformed quantities
\begin{eqnarray*}
&&\bar{u}(t,x)= u(t,\Xi_h(t,x)),\quad \bar{\pi}(t,x)=\pi(t,\Xi_h(t,x)),
\quad t>0,\; x\in\Omega\backslash\Sigma,\\
&&\bar{\theta}(t,x)=\theta(t,\Xi_h(t,x)),\quad t>0,\; x\in\Omega\backslash\Sigma,\\
&&\bar{j}(t,x)=j(t,\Xi_h(t,x)),\quad\ t>0,\; x\in\Sigma,
\end{eqnarray*}
the {\em pull back} of $(u,\pi,\theta)$ and $j$. This way we have transformed the time varying regions $\Omega\setminus\Gamma(t)$ to the fixed domain
$\Omega\setminus\Sigma$.

This transformation gives the following problem for $(\bar{u},\bar{\pi}, \bar{\theta}, h)$:
\begin{equation}
\label{ti2pp}
\begin{aligned}
\partial_t \bar{u} -\cG(h)\cdot\mu(\bar{\theta})(\cG(h)\bar{u}+[\cG(h)\bar{u}]^{\sf T}) +\cG(h)\bar{\pi}&=
\cR_u(\bar{u},\bar{\theta},h) &&\text{in}&&\Omega\backslash\Sigma,\\
\cG(h)\cdot\bar{u}&=0 &&\text{in}&&\Omega\backslash\Sigma,\\
\kappa(\bar{\theta})\partial_t \bar{\theta}-\cG(h)\cdot d(\bar{\theta})\cG(h)\bar{\theta}&=
\cR_\theta(\bar{u},\bar{\theta},h) && \text{in} &&\Omega\backslash\Sigma, \\
\bar{u}=\partial_\nu \bar{\theta}&=0 && \text{on}&&\partial\Omega,\\
-[\![\mu(\bar{\theta})(\cG(h)\bar{u}+[\cG(h)\bar{u}]^{\sf{T}})\nu_\Gamma(h)
-\bar{\pi}\nu_\Gamma(h)]\!] &=
\sigma H_\Gamma(h) \nu_\Gamma(h) &&\text{on} && \Sigma,\\
[\![\bar{u}]\!]=[\![\bar{\theta}]\!]&=0 && \text{on}&&\Sigma,\\
[\![\psi(\bar{\theta})]\!] + \sigma H_\Gamma(h)&=0 &&\text{on}&& \Sigma,\\
\beta(h)\partial_t h - \bar{u}\cdot \nu_\Gamma -[\![d(\bar{\theta})\cG(h)\bar{\theta}\cdot\nu_\Gamma]\!]/l(\bar{\theta})&=0 &&\text{on}&&\Sigma, \\
\bar{u}(0)=\bar{u}_0,\quad\bar{\theta}(0)=\bar{\theta}_0,\quad h(0)&=h_0.&&
\end{aligned}
\end{equation}
Here $\cG(h)$ and $H_\Gamma(h)$ denote the transformed gradient and curvature, respectively.
More precisely
we have
$$D\,\Xi_h  = I + D\xi_h, \quad \quad [D\,\Xi_h]^{-1}
= I - {[I + D\xi_h]}^{-1}D\xi_h=:I-M_1(h)^{\sl T}, $$
with
\begin{align*}
D\xi_h(t,x)&= \frac{1}{a}\chi\prime(d_\Sigma(x)/a)h(t,\Pi_\Sigma(x))\nu_\Sigma(\Pi_\Sigma(x))\otimes\nu_\Sigma(\Pi_\Sigma(x))\\
& + \chi(d_\Sigma(x)/a)\nu_\Sigma(\Pi_\Sigma(x))\otimes M_0(d_\Sigma(x))\nabla_\Sigma h(t,\Pi_\Sigma(x))\\
&-\chi(d_\Sigma(x)/a)h(t,\Pi_\Sigma(x))L_\Sigma(\Pi(x))M_0(d_\Sigma(x))P_\Sigma(\Pi_\Sigma(x))\\
&-h(t,\Pi(x))L_\Sigma(\Pi(x))M_0(d_\Sigma(x))P_\Sigma(\Pi_\Sigma(x))
\quad \mbox{ for }\; |d_\Sigma(x)|<2a/3,\\
D\xi_h(t,x)&=0
\quad\mbox{ for }\; |d_\Sigma(x)|>2a/3.
\end{align*}
In particular,
\begin{align*}
D\xi_h(t,x)
&= \nu_\Sigma(\Pi(x))\otimes M_0(d_\Sigma(x))\nabla_\Sigma h(t,\Pi(x))
\quad\mbox{for } |d_\Sigma(x)|<a/3.
\end{align*}
Thus, $[I + D\xi_h]$
is boundedly invertible if $h$ and $\nabla_\Sigma h$ are sufficiently small, e.g.\ if
\begin{equation}
\label{hanzawainv}
{|h|}_\infty < \frac{1}{3}\min\{a/|\chi\prime|_\infty,1/|L_\Sigma|_\infty\}
\quad \mbox{ and }\; {|\nabla_\Sigma h|}_\infty < \frac{1}{3}.
\end{equation}
With these properties we derive
\begin{align*}
\nabla\pi\circ\Xi_h
&= \cG(h)\bar{\pi}=  [\nabla\Xi_h^{-1}\circ\Xi_h]^{\sf T}\nabla\bar{\pi}
 = [\nabla\Xi_h]^{-1,\sf T}\nabla\bar{\pi}
 = (I - M_1(h))\nabla\bar{\pi}\\
 \text{div}\,u\circ\Xi_h
 &= \cG(h)\cdot\bar{u}=  (I - M_1(h))\nabla\cdot \bar{u}.
\end{align*}
Next we note that
\begin{align*}
\partial_t u\circ\Xi_h
&=  \partial_t\bar{u} -  [D u\circ\Xi_h]\partial_t\Xi_h
=  \partial_t\bar{u} -  D\bar u\,[D\,\Xi_h]^{-1}\,\partial_t\xi_h  \\
& =  \partial_t\bar{u} - ([I+ D\xi_h]^{-1}\,\partial_t\xi_h\cdot \nabla)u
 =: \partial_t\bar{u} + (R(h)\cdot \nabla)\bar{u}.
\end{align*}
with
$
R(h)=-[I+D\xi_h]^{-1}\,\partial_t\xi_h.
$
Hence
$$\cR_u(\bar{u},\bar{\theta},h)=-\bar{u}\cdot\cG(h)\bar{u}-(R(h)\cdot\nabla)\bar{u}.$$
Observe that the function $R(h)$ contains a time derivative of $h$ linearly.
Similarly we get
\begin{align*}
\nabla\theta\circ\Xi_h       &=\cG(h)\bar{\theta}=(I - M_1(h))\nabla\bar{\theta}\\
\partial_t \theta\circ \Xi_h &=
 \partial_t\bar\theta - \nabla\bar\theta\cdot [\nabla\Xi_h]^{-1}\partial_t\xi_h
=:\partial_t\bar{\theta} + (R(h)\cdot\nabla)\bar{\theta},
\end{align*}
and so
$$\cR_\theta(\bar{u},\bar{\theta},h))=-\kappa(\bar{\theta})\bar{u}\cdot\cG(h)\bar{\theta}-\kappa(\theta)
(R(h)\cdot\nabla)\bar{\theta}
+\mu(\bar{\theta})\big(\cG(h)\bar{u}+[\cG(h) \bar{u}]^T\big):\cG(h)\bar{u}.$$
With the Weingarten tensor $L_\Sigma$ and the surface gradient
 $\nabla_\Sigma$ we further have
\begin{eqnarray*}
\nu_\Gamma(h)= \beta(h)(\nu_\Sigma-\alpha(h)),&& \alpha(h)= M_0(h)\nabla_\Sigma h,\\
M_0(h)=(I-hL_\Sigma)^{-1},&& \beta(h) = (1+|\alpha(h)|^2)^{-1/2},
\end{eqnarray*}
and
$$V_\Gamma=\partial_t\Xi\cdot \nu_\Gamma = \partial_t h \nu_\Gamma\cdot \nu_\Sigma=
\beta(h)\partial_t h.$$
The curvature $H_\Gamma(h)$ becomes
$$ H_\Gamma(h) = \beta(h)\{ {\rm tr} [M_0(h)(L_\Sigma+\nabla_\Sigma \alpha(h))]
-\beta^2(h)M_0(h)\alpha(h)\cdot[\nabla_\Sigma\alpha(h)]\alpha(h)\},$$
a differential expression involving second order derivatives of $h$ only linearly.
Its linearization is given by
$$H^\prime_\Gamma(0)= {\rm tr}\, L_\Sigma^2 +\Delta_\Sigma.$$
Here $\Delta_\Sigma$ denotes the Laplace-Beltrami operator on $\Sigma$.

It is convenient to decompose the stress boundary condition into tangential and normal parts.
For this purpose let as before $P_{\Sigma}= I -\nu_{\Sigma}\otimes\nu_{\Sigma}$ denote the projection
onto the tangent space of $\Sigma$. Multiplying the stress interface condition with $\nu_\Sigma/\beta$ we obtain
$$[\![\bar{\pi}]\!]- \sigma H_\Gamma(h)=
[\![\mu(\theta)(\cG(h)\bar{u}+[\cG(h)\bar{u}]^{\sf T})(\nu_\Sigma-M_0(h)\nabla_\Sigma h)\cdot\nu_\Sigma ]\!],
$$
for the normal part of the stress interface condition.
Substituting this expression for $[\![\bar{\pi}]\!]- \sigma H_\Gamma(h)$ in \eqref{ti2pp}
and applying the tangential projection yields
\begin{align*}
&-P_\Sigma[\![\mu(\theta)(\cG(h)\bar{u}+[\cG(h)\bar{u}]^{\sf T})(\nu_\Sigma-M_0(h)\nabla_\Sigma h)]\!]\\
&\qquad = [\![\mu(\theta)(\cG(h)\bar{u}
+[\cG(h)\bar{u}]^{\sf T})(\nu_\Sigma-M_0(h)\nabla_\Sigma h)\cdot\nu_\Sigma ]\!] M_0(h)\nabla_\Sigma h
\end{align*}
for the tangential part. Note that the latter neither contains the pressure jump nor the curvature which is the advantage of this decomposition.

The idea of our approach can be described as follows.
We consider the transformed problem \eqref{ti2pp}.
Based on maximal $L_p$-regularity of the linear problem given by the left
hand side of \eqref{ti2pp}, we employ the contraction mapping principle to obtain local
well-posedness of the nonlinear problem.
This program will be carried out in the next sections.

%%%%%%%%%%%%%%%%%%%%%%%%%%%%%%%%
\section{The Linear Problem}
The principal part of the linearization of \eqref{ti2pp} reads as follows.
\begin{equation}
\label{princlin-u}
\begin{aligned}
\partial_t u-\mu_0(x) \Delta u +\nabla\pi &=f_u &&\text{in}&& \Omega\setminus\Sigma,\\
{\rm div}\, u &=f_d &&\text{in}&& \Omega\setminus\Sigma,\\
[\![u]\!] &=0 &&\text{on} &&\Sigma,\\
 -[\![\mu_0(x)(\nabla u+[\nabla u]^T)\nu_\Sigma-\pi\nu_\Sigma]\!] &= \sigma \Delta_\Sigma h\nu_\Sigma +g_u
&&\text{on}&& \Sigma,\\
u&=0 &&\text{on}&& \partial\Omega,\\
u(0)&=u_0 &&\text{in}&& \Omega,
\end{aligned}
\end{equation}
\vspace{1mm}
\begin{equation}
\label{princlin-theta}
\begin{aligned}
\kappa_0(x)\partial_t\vtheta -d_0(x)\Delta \vtheta &=f_\theta &&\text{in}&& \Omega\setminus\Sigma,\\
[\![\vtheta]\!]&=0 &&\text{on}&&\Sigma,\\
\partial_\nu\vtheta &=0 &&\text{on}&&\partial\Omega,\\
\vtheta(0)&=\vtheta_0 &&\text{in}&& \Omega,
\end{aligned}
\end{equation}
\vspace{1mm}
\begin{equation}
\label{princlin-Gamma}
\begin{aligned}
 l_1(t,x)\vtheta + \sigma\Delta_\Sigma h  &= g_\theta &&\text{on}&&\Sigma,\\
 \partial_t h - [\![d_0(x)\partial_\nu\vtheta]\!]/l_0(x) &=g_h &&\text{on}&&\Sigma,\\
 h(0)&=h_0 &&\text{on}&& \Sigma.
\end{aligned}
\end{equation}
Here
\begin{equation*}
\begin{split}
&\mu_0(x)=\mu(\theta_0(x)),\quad \kappa_0(x)=\kappa(\theta_0(x)),\quad d_0(x)=d(\theta_0(x)),\\
&l_0(x)=l(\theta_0(x)),\quad
l_1(t,x)=[\![\psi^\prime((e^{\Delta_\Sigma t}\gamma\theta_0)(x))]\!],
\end{split}
\end{equation*}
where $\gamma\theta_0$ means the restriction of $\theta_0$ to $\Sigma$. Observe that the term $u\cdot\nu_\Gamma$ in the equation for $h$ is of lower order as it enjoys more regularity than the trace of $\theta$ on $\Sigma$.
Therefore this system is triangular, since $u$ does neither appear
in \eqref{princlin-theta} nor in \eqref{princlin-Gamma}.
The latter system comprises the linearized Stefan problem with
surface tension which has been studied in
 \cite[Theorem 3.3]{PSZ10}.
To state this result we define the solution spaces
\begin{align*}
\EE_\theta(J)&= \{\vtheta\in H^1_p(J;L_p(\Omega))\cap
L_p(J;H^2_p(\Omega\setminus\Sigma)\cap C(\bar{\Omega})): \partial_{\nu}\vtheta=0\text{ on }\partial\Omega \},\\
\EE_h(J)&= W^{3/2-1/2p}_p(J;L_p(\Sigma))\cap
W^{1-1/2p}_p(J;H^2_p(\Sigma))\cap L_p(J;W^{4-1/p}_p(\Sigma)),\\
\EE(J)&=\EE_\theta(J)\times\EE_h(J)
\end{align*}
and the spaces of data
\begin{align*}
\FF_\theta(J)&=L_p(J\times\Omega),\\
\FF_H(J)&=W^{1-1/2p}_p(J;L_p(\Sigma))\cap L_p(J;W^{2-1/p}_p(\Sigma)),\\
\FF_h(J)&=W^{1/2-1/2p}_p(J;L_p(\Sigma))\cap L_p(J;W^{1-1/p}_p(\Sigma)),\\
\FF(J)&=\FF_\theta(J)\times\FF_H(J)\times\FF_h(J).
\end{align*}
Then we have
%%%%%%%%%%%%%%%%%%%%%%%%
\begin{theorem}
\label{linpeq0} Let $p>n+2$ and $\sigma>0$, suppose
$\kappa_0\in C(\bar{\Omega}_i)$, $d_0\in C^1(\bar{\Omega}_i)$, $i=1,2$, $\kappa_0,d_0>0$
on $\bar{\Omega}$, $l_0\in C^1(\Sigma)$, and let
$l_1\in \FF_H(J)$
such that $l_0l_1>0$ on $J\times\Sigma$, where $J=[0,t_0]$ is a
finite time interval.
\\
Then there is a unique solution $(\vtheta,h)\in\EE(J)$
of \eqref{princlin-theta}-\eqref{princlin-Gamma} if and only if
the data $(f_\theta,g_\theta,g_h,\vtheta_0,h_0)$ satisfy
$$(f_\theta,g_\theta,g_h)\in \FF(J),\quad
(\vtheta_0,h_0)\in [W^{2-2/p}_p(\Omega\setminus\Sigma)\cap C(\bar{\Omega})]
\times W^{4-3/p}_p(\Sigma),$$
and the compatibility conditions
\begin{equation*}
\partial_{\nu}\vtheta_0=0\text{ on }\partial\Omega,\quad
l_1(0)\vtheta_0+\sigma \Delta_\Sigma h_0=g_\theta(0),
\quad g_h(0) +[\![d_0\partial_\nu \vtheta_0]\!]/l_0
\in W^{2-6/p}_p(\Sigma).
\end{equation*}
The solution map
$[(f_\theta,g_\theta,g_h,z_0)\mapsto (\vtheta,h)]$
is continuous between the corresponding spaces.
\end{theorem}
%%%%%%%%%%%%%%%%%%%%%%%
Having solved this problem, $\Delta_\Sigma h\in \FF_h(J)$,
hence we may now solve the remaining Stokes problem to obtain
a unique solution $(u,\pi)$ in the class
$$ u\in H^1_p(J;L_p(\Omega))^n\cap L_p(J;H^2_p(\Omega\setminus\Sigma)
\cap C(\bar{\Omega}))^n,\quad \nabla\pi\in L_p(J,L_p(\Omega))^n,$$
provided the data satisfy
$$f_u\in \FF_\theta(J)^n,\quad g_u\in \FF_h(J)^n,\quad
f_d\in H^1_p(J;\dot{H}^{-1}_p(\Omega))\cap
L_p(J;H^1_p(\Omega\setminus\Sigma)),$$
as well as $u_0\in [W^{2-2/p}_p(\Omega\setminus\Sigma)\cap C(\bar{\Omega})]^n$,
and the compatibility conditions
\begin{equation*}
u_0=0\text{ on } \partial\Omega,
\quad{\rm div}\, u_0=f_d(0),
\quad -P_\Sigma[\![\mu_0(\nabla u_0+[\nabla u_0]^{\sf T})\nu_\Sigma]\!]
= P_\Sigma g_u(0),
\end{equation*}
are satisfied. Concerning the two-phase Stokes problem we refer to
Shibata and Shimizu \cite{ShSh09} and to K\"ohne, Pr\"uss and Wilke
\cite{KPW10},
where the case $\mu_0=const$ has been treated.
This result can be
extended to nonconstant $\mu_0$ via the method of localization.
For the one-phase case this has been carried out for much more
general Stokes problems in Bothe and Pr\"uss \cite{BoPr07}

Therefore the linearized problem \eqref{princlin-u}-\eqref{princlin-Gamma} has the property of
maximal $L_p$-regularity. To state this result we set
\begin{align*}
&\EE_u(J):=\{u\in [H^1_p(J;L_p(\Omega))\cap
L_p(J;H^2_p(\Omega\setminus\Sigma)\cap C(\bar{\Omega}))]^n: u=0\text{ on } \partial\Omega \},\\
&\EE_\pi(J):= L_p(J;\dot{H}_p^1(\Omega\setminus\Sigma)),
\end{align*}
and define the solution space for \eqref{princlin-u}-\eqref{princlin-Gamma} as
$$\EE(J)= \EE_u(J)\times\EE_\pi(J)\times \FF_h(J)\times\EE_\theta(J)\times\EE_h(J).$$
Then the main result on the linearized problem reads
%%%%%%%%%%%%%%%%%%%%%
\begin{theorem}
\label{maximalfull}
Let $p>n+2$ and $\sigma>0$.
Suppose $\mu_0$, $\kappa_0\in C(\bar\Omega_i)$, $d_0\in C^1(\bar\Omega_i)$, $i=1,2$,
$\kappa_0,d_0>0$ on $\bar{\Omega}$, $l_0\in C^1(\Sigma)$, and
$l_1\in \FF_h(J)$
such that $l_0l_1>0$ on $J\times\Sigma$, where $J=[0,t_0]$ is a
finite time interval.

Then the linear problem \eqref{princlin-u}-\eqref{princlin-Gamma} admits a unique solution
$(u,\pi,[\![\pi]\!],\theta, h)\in\EE(J)$
if and only if the data $(u_0,\vtheta_0,h_0)$ and $(f_u,f_d,g_u,f_\theta,
g_\theta,g_h)$ satisfy the regularity conditions:
\begin{align*}
&(u_0,\vtheta_0,h_0)\in W_p^{2-2/p}(\Omega\setminus\Sigma)^n\times
W_p^{2-2/p}(\Omega\setminus\Sigma)\times W_p^{4-3/p}(\Sigma),\\
&(f_u,f_\theta) \in \FF_\theta(J)^{n+1},\quad
f_d \in H_p^1(J; \dot H_p^{-1}(\Omega)) \cap L_p(J; H_p^1(\Omega
\setminus \Sigma)),\\
&(g_u,g_h)\in \FF_h(J)^{n+1},
\quad g_\theta\in \FF_H(J),
\end{align*}
and the compatibility conditions:
\begin{equation*}
\begin{aligned}
{\rm div}\, u_0&=f_d(0) &&{\rm in}&& \Omega\setminus\Sigma,\\
u_0=\pd_\nu\theta_0&=0 &&{\rm on}&& \pd\Omega,\\
[\![u_0]\!]=[\![\theta_0]\!]&=0 &&{\rm on}&& \Sigma,\\
-P_\Sigma[\![\mu_0(\nabla u_0+[\nabla u_0]^{\sf T})\nu_\Sigma]\!]&=
P_\Sigma g_u(0) &&{\rm on}&& \Sigma,\\
l_1(0)\vtheta_0+\sigma\Delta_\Sigma h_0&=g_\theta(0) &&{\rm on}&&\Sigma\\
 g_h(0)+ [\![d_0 \pd_\nu\theta_0]\!]/l_0 &\in W_p^{2-6/p}(\Sigma). &&
\end{aligned}
\end{equation*}
The solution map
$[(u_0, \theta_0, h_0, f_u, f_d, g_u, f_\theta, g_\theta, g_h)
\mapsto (u, \pi, [\![\pi]\!], \vtheta, h)]$ is  continuous between the corresponding spaces.
\end{theorem}

\section{Local Well-Posedness}

The basic result for local well-posedness of problem \eqref{i2pp}
in an $L_p$-setting is the following.
%%%%%%%%%%%%%%%%%%%%%
\begin{theorem} \label{wellposed} Let $p>n+2$,  $\sigma>0$, suppose
$\psi_i\in C^3(0,\infty)$, $\mu_i,d_i\in C^2(0,\infty)$ such that
$$\kappa_i(s)=-s\psi_i^{\prime\prime}(s)>0,\quad \mu_i(s)>0,
\quad  d_i(s)>0,\quad s\in(0,\infty),\; i=1,2.$$
Let $\Omega\subset\R^n$ be a bounded domain with boundary
$\partial\Omega\in C^{3-}$ and suppose $\Gamma_0\subset \Omega$ is a
closed hypersurface.
Assume the {\bf regularity conditions}
\begin{equation*}
u_0\in W^{2-2/p}_p(\Omega\setminus\Gamma_0)^n,\quad
\theta_0\in W^{2-2/p}_p(\Omega\setminus\Gamma_0),\quad
\Gamma_0\in W^{4-3/p}_p,
\end{equation*}
the {\bf  compatibility conditions}
\begin{equation*}
\begin{aligned}
{\rm div}\, u_0&=0 &&{\rm in} &&\Omega\setminus\Gamma_0,\\
u_0=\partial_\nu\theta_0 &=0 &&{\rm on}&&\partial\Omega,\\
[\![u_0]\!]=[\![\theta_0]\!]&=0 &&{\rm on}&&\Gamma_0,\\
P_{\Gamma_0}
[\![\mu_0(\nabla u_0+[\nabla u_0]^{\sf T})\nu_{\Gamma_0}]\!] &=0 &&{\rm on}&&\Gamma_0,\\
[\![\psi(\theta_0)]\!]+\sigma H_{\Gamma_0}&=0 &&{\rm on}&& \Gamma_0,\\
[\![d_0\partial_\nu \theta_0]\!]\in
W^{2-6/p}_p&(\Gamma_0), &&
\end{aligned}
\end{equation*}
and the {\bf well-posedness condition}
$$ l(\theta_0)\neq0\; \mbox{ on }\; \Gamma_0\quad
\mbox{ and }\;\theta_0>0 \; \mbox{ on } \;\bar{\Omega}.$$
Then there exists a unique $L_p$-solution of problem \eqref{i2pp}  on
some possibly small but nontrivial time interval $J=[0,a]$.
\end{theorem}
%%%%%%%%%%%%%%%%%%%%%%

\medskip

\noindent
Here the notation $\Gamma_0\in W^{4-3/p}_p$ means that $\Gamma_0$ is
a $C^2$-manifold, such that
its (outer) normal field $\nu_{\Gamma_0}$ is of class
$W^{3-3/p}_p(\Gamma_0)$. Therefore the Weingarten tensor
 $L_{\Gamma_0}=-\nabla_{\Gamma_0}\nu_{\Gamma_0}$ of $\Gamma_0$
 belongs to $W^{2-3/p}_p(\Gamma_0)$ which embeds into
$C^{1+\alpha}(\Gamma_0)$, with $\alpha=1-(n+2)/p>0$ since $p>n+2$
by assumption.
For the same reason we also have
$u_0\in C^{1+\alpha}(\bar{\Omega}_i(0)))^n$, and
$\theta_0\in C^{1+\alpha}(\bar{\Omega}_i(0)))$, $i=1,2$,
and $V_0\in C^{2\alpha}(\Gamma_0)$. The notion $L_p$-solution means
that $(u,\pi,\theta,\Gamma)$ is obtained as the push-forward of an
$L_p$-solution of the transformed problem \eqref{ti2pp}.
The proof of Theorem \ref{wellposed} is given in the next section.

\bigskip

\noindent
For later use we discuss an extension of the local existence results
to spaces with time weights. For this purpose, given a
$U\!M\!D$-Banach space $Y$
and $\mu\in(1/p,1]$, we define for $J=(0,t_0)$
$$K^s_{p,\mu}(J;Y):=\{u\in L_{p,loc}(J;Y):
\; t^{1-\mu}u\in K^s_p(J;Y)\},$$
where $s\geq0$ and $K\in\{H,W\}$. It has been shown in
Pr\"uss and Simonett \cite{PrSi04} that the operator $d/dt$ in
$L_{p,\mu}(J;Y)$ with domain
$$D(d/dt)={_0H}^1_{p,\mu}(J;Y)=\{u\in H^1_{p,\mu}(J;Y):\; u(0)=0\}$$
is sectorial and admits an $H^\infty$-calculus with angle $\pi/2$.
This is the main tool to extend the results for the linear problem,
i.e.\  Theorem \ref{maximalfull}, to the time weighted setting,
where the solution space $\EE(J)$ is replaced by $\EE_{\mu}(J)$
and $\FF(J)$ by $\FF_\mu(J)$ respectively, where
\begin{equation*}
z\in \EE_{\mu}(J)\Leftrightarrow t^{1-\mu}z\in \EE(J),\quad
f\in \FF_{\mu}(J)\Leftrightarrow t^{1-\mu}f\in \FF(J).
\end{equation*}
The trace spaces for $(u,\theta)$  and $h$ for $p>3$ are then given by
\begin{equation}
\label{tracesp-mu}
\begin{aligned}
&(u_0,\theta_0)\in [W^{2\mu-2/p}_p(\Omega\setminus\Sigma)\cap C(\bar\Omega)]^{n+1},
\quad  h_0\in W^{2+2\mu-3/p}_p(\Sigma), \\
& h_1:=(\partial_th)(0)\in W^{4\mu-2-6/p}_p(\Sigma),
\end{aligned}
\end{equation}
where for the last trace  we need in addition $\mu>1/2+3/2p$.
Note that the embeddings
\begin{equation*}
\begin{aligned}
\EE_{\mu,u}(J)\times \EE_{\mu,\theta}(J)\hookrightarrow [C(J\times \bar{\Omega})
\cap C(J;C^1(\bar{\Omega}_i)]^{n+1},
\quad \quad\EE_{\mu,h}(J)\hookrightarrow C(J;C^3(\Sigma))
\end{aligned}
\end{equation*}
require
$\mu>1/2+(n+2)/2p$, which is feasible since $p>n+2$ by assumption.
This restriction is needed for the estimation of the nonlinearities.

The assertions for the linear problem remain valid for such $\mu$,
replacing $\EE(J)$ by $\EE_\mu(J)$, $\FF(J)$ by $\FF_\mu(J)$,
for initial data subject to \eqref{tracesp-mu}.
This relies on the fact mentioned above that $d/dt$ admits a bounded
$H^\infty$-calculus
with angle $\pi/2$ in the spaces $L_{p,\mu}(J;Y)$.
Concerning such time weights, we refer to Pr\"uss, Simonett and
Zacher \cite{PSZ10} for the Stefan problem, and to K\"ohne, Pr\"uss
and Wilke \cite{KPW10} for the isothermal incompressible two-phase
Navier-Stokes problem, and also to Meyries and Schnaubelt \cite{MeSc11} for a general theory.
Thus as a consequence of these considerations we have the following result.
%%%%%%%%%%%%%%%%%%%%%%
\begin{corollary}
\label{wellposed3}
Let $p>n+2$,
$\mu\in (1/2+(n+2)/2p,1]$, $\sigma>0$, and suppose
$\psi\in C^3(0,\infty)$, $\mu,d\in C^2(0,\infty)$ are such that
$$\kappa_i(s)=-s\psi_i^{\prime\prime}(s)>0,\quad \mu_i(s)>0,\quad
d_i(s)>0,\quad s\in(0,\infty),\;i=1,2.$$
Let $\Omega\subset\R^n$ be a bounded domain with boundary
$\partial\Omega\in C^{3-}$ and suppose $\Gamma_0\subset \Omega$
is a closed hypersurface.
Assume the {\bf regularity conditions}
$$(u_0,\theta_0)\in [W^{2\mu-2/p}_p(\Omega\setminus\Gamma_0)\cap
C(\bar{\Omega})]^{n+1},\quad \Gamma_0\in W^{2+2\mu-3/p}_p,$$
 the {\bf  compatibility conditions}
\begin{equation*}
\begin{aligned}
{\rm div}\, u_0&=0 &&{\rm in} &&\Omega\setminus\Gamma_0,\\
u_0=\partial_\nu\theta_0 &=0 &&{\rm in}&&\partial\Omega,\\
[\![u_0]\!]=[\![\theta_0]\!]=0,\quad P_{\Gamma_0}
[\![\mu_0(\nabla u_0+[\nabla u_0]^{\sf T})\nu_{\Gamma_0}]\!] &=0 &&{\rm on}&&\Gamma_0,\\
[\![\psi(\theta_0)]\!]+\sigma H_{\Gamma_0}&=0 &&{\rm on}&& \Gamma_0,\\
\quad [\![d_0\partial_{\nu} \theta_0]\!]\in
W^{4\mu-2-6/p}_p&(\Gamma_0), &&
\end{aligned}
\end{equation*}
as well as  the {\bf well-posedness condition}
$$ l(\theta_0)\neq0\; \mbox{ on }\; \Gamma_0\quad
\mbox{ and }\;\theta_0>0 \; \mbox{ on } \;\bar{\Omega}.$$
Then the transformed problem \eqref{ti2pp} admits a unique solution
$(u,\pi,[\![\pi]\!],\vtheta, h)\in\EE_\mu((0,a))$ for some
$a>0$. The solution depends continuously on the data.
For each $\delta\in(0,a)$, the solution belongs to $\EE(\delta,a)$,
i.e.\ regularizes instantly.
\end{corollary}
%%%%%%%%%%%%%%%%%%%%%%%%%
\medskip
%\goodbreak

\section{Proof of the Main Result}

In this section we prove Theorem \ref{wellposed}, for given initial data
$\Gamma_0\in W_p^{4-3/p}$, $u_0\in W_p^{2-2/p}
(\Omega\setminus\Gamma_0)^n$
and $\theta_0\in W_p^{2-2/p}(\Omega\setminus\Gamma_0)$
satisfying the compatibility conditions and the well-posedness condition
stated in Theorem \ref{wellposed}. We suppose that
$\psi_i\in C^3(0,\infty)$, $\mu_i,d_i\in C^2(0,\infty)$
satisfy
$$
\kappa_i(s)=-s\psi_i^{\prime\prime}(s)>0,\quad \mu_i(s)>0,\quad
d_i(s)>0,\quad s\in(0,\infty),\; i=1,2.
$$
According to the considerations in Section 4, $\Gamma_0$ can be approximated
by a real analytic hypersurface for any prescribed $\eta>0$ in the sense
that $d_H(\cN^2\Sigma, \cN^2\Gamma_0)<\eta$, and is parameterized by
$h_0\in W_p^{4-3/p}(\Sigma)$. It is sufficient to prove the local
well-posedness of the nonlinear problem \eqref{ti2pp}, because the
notion of $L_p$-solution means that the solution of \eqref{i2pp} is
obtained as the push-forward of an $L_p$-solution of the transformed
problem \eqref{ti2pp}.
In order to facilitate this task, we rewrite \eqref{ti2pp} in
quasilinear form, dropping the bars and collecting its principal linear
part on the left hand side. We set as before $\mu_0(x)=\mu(\theta_0(x))$, $\kappa_0(x)=\kappa(\theta_0(x))$, $d_0(x)=d(\theta_0(x))$, $l_0(x)=l(\theta_0(x))$,
and $l_1(t,x)= [\![\psi^\prime(e^{\Delta_\Sigma t}\gamma\theta_0(x))]\!]$.
We then have

%{\allowdisplaybreaks
\begin{equation}
\label{quasi}
\begin{aligned}
\partial_t u-\mu_0(x) \Delta u +\nabla\pi
&=F_u(u,\pi,\theta,h) &&\text{in}&& \Omega\setminus\Sigma,\\
{\rm div}\, u
&=F_d(u,h)&&\text{in}&& \Omega\setminus\Sigma,\\
[\![u]\!] &= 0 &&\text{on}&& \Sigma,\\
-P_\Sigma[\![\mu_0(x)(\nabla u+[\nabla u]^{\sf T}) \nu_\Sigma]\!]
&= G_u^{tan}(u,\theta,h) &&\text{on}&&\Sigma,\\
-2[\![\mu_0(x)\nabla u]\!] \nu_\Sigma\cdot \nu_\Sigma
+[\![\pi]\!]- \sigma \Delta_\Sigma h
&= G_u^{nor}(u,\theta,h) &&\text{on}&& \Sigma,\\
u=\partial_\nu\theta &=0 &&\text{on}&&\partial\Omega,\\
\kappa_0(x)\partial_t\theta -d_0(x)\Delta \theta
&=F_\theta(u,\theta,h) &&\text{on}&&\Omega\setminus\Sigma,\\
[\![\theta]\!]&=0 &&\text{on}&& \Sigma,\\
 l_1(t,x)\theta +\sigma\Delta_\Sigma h
&= G_\theta(\theta,h) &&\text{on}&&\Sigma,\\
\partial_t h  -[\![d_0(x)\partial_\nu\theta]\!]/l_0(x) &=G_h(u,\theta,h) &&\text{on}&& \Sigma,\\
u(0)=u_0,\quad\theta(0)&=\theta_0 && \text{in}&& \Omega,\\
h(0)&=h_0 &&\text{in} && \Sigma.
\end{aligned}
\end{equation}
%}
The nonlinearities are defined by
{\allowdisplaybreaks
\begin{align*}
F_u(u,\theta,\pi,h)
&=(\mu(\theta)\!-\!\mu(\theta_0))\Delta u + M_1(h)\nabla\pi - (u\cdot (I-M_1(h))\nabla)u\!-\!(R(h)\cdot\nabla)u\\
&+\mu^\prime(\theta)
\big( (I-M_1(h)) \nabla\theta\cdot ((I-M_1(h))\nabla u +[(I-M_1(h))\nabla u]^{\sf T})\big) \\
&-\mu(\theta)(M_2(h):\nabla^2) u-\mu(\theta)(M_3(h)\cdot\nabla) u+\mu(\theta)M_4(h):\nabla u,\\
F_d(u,h)&=M_1(h):\nabla u,\\
%\end{align*}
%\begin{align*}
G_u^{tan}(u,h)
&= P_\Sigma[\![(\mu(\theta)-\mu(\theta_0))(\nabla u +[\nabla u]^T)\nu_\Sigma]\!]\\
&-P_\Sigma [\![\mu(\theta)(\nabla u+[\nabla u]^{\sf T})
M_0(h)\nabla_\Sigma h]\!]\\
&-P_\Sigma [\![\mu(\theta)(M_1(h)\nabla u+[M_1(h)\nabla u]^{\sf T})
(\nu_\Sigma-M_0(h)\nabla_\Sigma h)]\!]\\
&+[\![\mu((I-M_1)\nabla u+[(I-M_1)\nabla u]^T)
(\nu_\Sigma-M_0\nabla_\Sigma h)\cdot\nu_\Sigma]\!]M(h)\nabla_\Sigma h,
\\
G_u^{nor}(u,h)&=
[\![(\mu(\theta)-\mu(\theta_0))(\nabla u +[\nabla u]^T)\nu_\Sigma\cdot\nu_\Sigma ]\!]\\
&-[\![\mu(\theta)(\nabla u+[\nabla u]^{\sf T})
M_0(h)\nabla_\Sigma h \cdot \nu_\Sigma ]\!]\\
&-[\![\mu(\theta)(M_1(h)\nabla u+[M_1(h)\nabla u]^{\sf T})
 (\nu_\Sigma-M_0(h)\nabla_\Sigma h )\cdot \nu_\Sigma ]\!]\\
&+\sigma(H_\Gamma(h)-\Delta_\Sigma h),\\
F_\theta(u,\theta,h)&= (\kappa(\theta_0)-\kappa(\theta))\partial_t \theta
-(d(\theta)-d(\theta_0))\Delta \theta
- d(\theta)M_2(h):\nabla^2 \theta\nonumber\\
& +d^\prime(\theta)|(I-M_1(h))\nabla \theta|^2 - d(\theta) M_3(h)\cdot\nabla \theta
-\kappa(\theta)(R(h)\cdot\nabla)\theta\nonumber\\
& -\kappa(\theta)u\cdot (I-M_1(h))\nabla\theta \\
& + \mu(\theta)((I-M_1(h))\nabla u+[(I-M_1(h))\nabla u]^T):(I-M_1(h))\nabla u,\\
G_\theta(\theta,h)&=l_1\theta-[\![\psi(\theta)]\!]
-\sigma (H_\Gamma(h)- \Delta_\Sigma h),
\label{FGdef}\\
G_h(u,\theta,h)&=[\![\big(d(\theta)/l(\theta)-d(\theta_0)/l(\theta_0)\big)\partial_\nu \theta]\!]
+u\cdot(\nu_\Sigma-M_0(h)\nabla_\Sigma h)\\
&-[\![(d/l) M_1(h)\nabla \theta\cdot(\nu_\Sigma-M_0(h)\nabla_\Sigma h)]\!]
-[\![(d/l)\nabla\theta\cdot M_0(h)\nabla_\Sigma h]\!].\nonumber
\end{align*}}
Here we employed the abbreviations
\begin{align*}M_2(h)&=M_1(h)+M_1^{\sf T}(h)-M_1(h)M_1^{\sf T}(h),\\
M_3(h)&=(I-M_1(h)):\nabla M_1(h), \\
M_4(h)&= ((I-M_1(h))\nabla) M_1(h)-[((I-M_1(h))\nabla) M_1(h)]^{\sf T}.
\end{align*}
We prove the local well-posedness of \eqref{quasi} by means of  maximal
$L_p$-regularity of the linear problem (Theorem \ref{maximalfull})
and the contraction mapping principle.
The right hand side of problem \eqref{quasi} consist of either lower order
terms, or terms of the same order as those appearing on the left hand side
but carry factors which can be made small by construction. Indeed, we have smallness of $h_0$,
$\nabla_\Sigma h_0$ and even of $\nabla^2_\Sigma h_0$ uniformly on
$\Sigma$, because $\Gamma_0$ is approximated by $\Sigma$ in the
second order bundle. $\theta$ appears nonlinearly in $\psi,\kappa,\mu,d$, but only to order zero; hence e.g.\ the difference
$\mu(\theta(t))-\mu(\theta_0)$ will be uniformly small for small times.
\par

We introduce appropriate function spaces. Let $J=[0,a]$.
The solution spaces are defined by
{\allowdisplaybreaks
\begin{align*}
&\EE_1(a):= \{u\in H_p^1(J; L_p(\Omega))^n\cap
L_p(J; H_p^2(\Omega\setminus \Sigma))^n :\, u=0\ {\rm on}\ \pd\Omega,
\ \ [\![u]\!]=0\},\\
&\EE_2(a):= L_p(J; \dot H^1_p(\Omega\setminus \Sigma)),\\
&\EE_3(a):= W_p^{1/2-1/2p}(J; L_p(\Sigma))
\cap L_p(J; W^{1-1/p}_p(\Sigma)),\\
&\EE_4(a):= \{\theta\in H_p^1(J; L_p(\Omega))\cap
L_p(J; H^2_p(\Omega\setminus \Sigma)):\, \pd_\nu\theta=0\
{\rm on}\ \pd\Omega,\ \ [\![\theta]\!]=0\},\\
&\EE_5(a):= W_p^{3/2-1/2p}(J; L_p(\Sigma))
\cap W_p^{1-1/2p}(J; H^2_p(\Sigma))
\cap L_p(J; W^{4-1/p}_p(\Sigma)).
\end{align*}}
We abbreviate
$$
\EE(a):= \{(u,\pi,[\![\pi]\!],\theta,h)\in
\EE_1(a) \times \EE_2(a) \times \EE_3(a) \times \EE_4(a) \times
\EE_5(a)\},
$$
and equip $\EE_j(a)$ ($j=1,\dots,5$) with their natural norms,
which turn $\EE(a)$ into a Banach space. A left subscript 0 always
means that the time trace at $t=0$ of the function in question is zero whenever it exists.
\par

The data spaces are defined by
{\allowdisplaybreaks
\begin{align*}
&\FF_1(a):= L_p(J; L_p(\Omega))^n,\\
&\FF_2(a):= H_p^{1}(J; \dot H_p^{-1}(\Omega)) \cap
L_p(J; H_p^1(\Omega\setminus\Sigma)),\\
&\FF_3(a):= W_p^{1/2-1/2p}(J; L_p(\Sigma))^n
\cap L_p(J; W^{1-1/p}_p(\Sigma))^n,\\
&\FF_4(a):= L_p(J; L_p(\Omega)),\\
&\FF_5(a):= W_p^{1-1/2p}(J; L_p(\Sigma))
\cap L_p(J; W^{2-1/p}_p(\Sigma)),\\
&\FF_6(a):= W_p^{1/2-1/2p}(J; L_p(\Sigma))
\cap L_p(J; W^{1-1/p}_p(\Sigma)).
\end{align*}}
We abbreviate
\begin{equation*}
\FF(a):= \{(f_u,f_d,g_u,g_\theta, g_\theta,g_h)\in \prod_{j=1}^6 \FF_j(a)\},
\end{equation*}
and equip $\FF_j(a)$ ($j=1,\dots,6$) with their natural norms,
which turn $\FF(a)$ into a Banach space.
\par
\medskip\noindent
{\bf Step 1.}\quad
In order to economize our notation,
we set
$z=(u,\pi,[\![\pi]\!], \theta,h)\in \EE(a)$
and reformulate the quasilinear problem \eqref{quasi} as
\begin{equation}
Lz=N(z)\quad (u(0),\theta(0),h(0))=(u_0,\theta_0,h_0),
\label{simple}
\end{equation}
where $L$ denotes the linear operator on the left hand side of
\eqref{quasi}, and $N$ denotes the nonlinear mapping on the right-hand
side of \eqref{quasi}. From Section~5 we know that $L:\EE(a)\to
\FF(a)$ is bounded and linear, and that $L:{}_0\EE(a)\to {}_0\FF(a)$
is an isomorphism for each $a>0$, with norm independent of $0<a\leq a_0<\infty$.
\par
Concerning the nonlinearity $N$, we have the following result.
%%%%%%%%%%%%%%%%%%%%%%%%%
\begin{proposition}
\label{pr:7.1}
Suppose $p>n+2$, $\sigma>0$, and let
$\psi_i\in C^3(0,\infty)$, $\mu_i,d_i\in C^2(0,\infty)$ such that
$$\kappa_i(s)=-s\psi_i^{\prime\prime}(s)>0,\quad \mu_i(s)>0,
\quad  d_i(s)>0,\quad s\in(0,\infty),\; i=1,2.$$
Then for each $a>0$ the nonlinearity satisfies
$N\in C^{1}(\EE(a),\FF(a))$
and its Fr{\'e}chet derivative $N^\prime$ satisfies in addition
$
N^\prime(u,\pi,[\![\pi]\!],\theta,h)\in \cL({}_0\EE(a),{}_0\FF(a)).
$
Moreover, there is $\eta>0$ such that
for a given $z_*\in\EE(a_0)$ with $|h_0|_{C^2(\Sigma)}\leq
\eta$, there are continuous functions $\alpha(r)>0$ and
$\beta(a)>0$ with $\alpha(0)=\beta(0)=0$,  such that
\begin{equation*}
 \Ver N^\prime(\bar{z}+z_*)\Ver_{\cB({_0\EE}(a),{_0\FF}(a))}
 \leq \alpha(r)+\beta(a), \quad \bar{z}\in \BB_r\subset {_0\EE}(a) .
\end{equation*}
\end{proposition}
\begin{proof}
The nonlinear terms $F_u(u,\theta,\pi, h)$, $F_d(u,h)$, $G_u^{tan}(u,\theta,h)$ and
$G_u^{nor}(u,\theta,h)$ are essentially the same as those for the isothermal two-phase
Navier-Stokes problem (cf.\ Proposition 4.2 in \cite{KPW10}).
The nonlinearities $F_\theta(u,\theta,h)$, $G_\theta(\theta,h)$, and
$G_h(u,\theta,h)$ are similar to those for $u\equiv0$ which have been discussed in \cite{PSZ10}. In fact, all appearing extra terms are of lower order.
\end{proof}
\par
\medskip\noindent
{\bf Step 2.}
We reduce the problem to initial values 0 and resolve the compatibilities as follows.
Thanks to Proposition 4.1 in \cite{KPW10}, we find
extensions
$f_d^*\in \FF_2(a)$, $g_u^*\in \FF_3(a)$
which satisfy
$$
f_d^*(0)={\rm div}\, u_0,\quad
P_\Sigma g_u^* = -P_\Sigma[\![ \mu_0(\nabla u_0+[\nabla u_0]^{\sf T}) \nu_\Sigma]\!].
$$
Then we define
$$g_\theta^*:= e^{\Delta_\Sigma t} [G_\theta(\theta,h)(0)],
\qquad g_h^*:= e^{\Delta_\Sigma t} G_h(u_0,\theta_0,h_0),$$
and  set $f_u^*=f_\theta^*=0$.
With these extensions, by Theorem 5.2 we may solve the linear problem
\eqref{princlin-u}-\eqref{princlin-Gamma} with initial data $(u_0,\theta_0,h_0)$ and inhomogeneities
$(f_u^*, f_d^*, g_u^*, f_\theta^*, g_\theta^*,g_h^*)$, which satisfy the required regularity conditions
and, by construction, the compatibility conditions,
to obtain a unique solution
$$z^*=(u^*,\pi^*,[\![\pi^*]\!],\theta^*,h^*)\in \EE(J)$$
with $u^*(0)=u_0$, $\theta^*(0)=\theta_0$, and $h^*(0)=h_0$.
\par
\medskip\noindent
{\bf Step 3.}\quad We rewrite problem \eqref{simple} as
$$
Lz=N(z+z^*)-Lz^*=: K(z), \quad z\in{}_0\EE(a).
$$
The solution is given by the fixed point problem $z=L^{-1}K(z)$, since Theorem~\ref{maximalfull}
implies that $L: {}_0\EE(a) \to {}_0\FF(a)$ is
an isomorphism with
$$
|L^{-1}|_{\cL({}_0\FF(a),{}_0\EE(a))} \le M, \quad a\in (0,a_0],
$$
where $M$ is independent of $a\le a_0$. We may assume that $M\ge 1$.
Thanks to Proposition
\ref{pr:7.1} and due to $K(0)=N(z^*)-Lz^*$, we may choose $a\in(0,a_0]$
and $r>0$ sufficiently small such that
$$
|K(0)|_{\FF(a)}\le \frac r{2M}, \quad
|K^\prime(z)|_{\cL({}_0\EE(a),{}_0\FF(a))} \le \frac 1{2M}, \quad
z\in {}_0\EE(a), \quad |z|_{\EE(a)}\le r
$$
hence
$$
|K(z)|_{\FF(a)}\le \frac rM,
$$
which ensures that $L^{-1}K(z): \BB_{{}_0\EE(a)}(0,r) \to
\BB_{{}_0\EE(a)}(0,r)$ is a contraction. Thus we may employ the contraction
mapping principle to obtain a unique solution on a time interval
$[0,a]$, which completes the proof of Theorem \ref{wellposed}.

\section{The Local Semiflow}
We follow here the approach introduced in K\"ohne, Pr\"uss and Wilke \cite{KPW10} for the isothermal incompressible two-phase Navier-Stokes problem without phase transitions and in Pr\"uss, Simonett and  Zacher \cite{PSZ10} for the Stefan problem with surface tension.

Recall that the closed $C^2$-hypersurfaces contained in $\Omega$ form a $C^2$-manifold,
which we denote by $\cMH^2(\Omega)$.
The charts are the parameterizations over a given hypersurface $\Sigma$ according to Section 4, and the tangent
space $T_\Sigma\cMH^2(\Omega)$ consists of the normal vector fields on $\Sigma$.
We define a metric on $\cMH^2(\Omega)$ by means of
$$d_{\cMH^2}(\Sigma_1,\Sigma_2):= d_H(\cN^2\Sigma_1,\cN^2\Sigma_2),$$
where $d_H$ denotes the Hausdorff metric on the compact subsets of $\R^n$. %% introduced in Section 3.
This way $\cMH^2(\Omega)$ becomes a Banach manifold of class $C^2$.

Let $d_\Sigma(x)$ denote the signed distance for $\Sigma$ as in Section 4.
We may then define the {\em level function} $\varphi_\Sigma$ by means of
$$\varphi_\Sigma(x) = \phi(d_\Sigma(x)),\quad x\in\R^n,$$
where
$$\phi(s)= s\chi(s/a) + ({\rm sgn}\,s)(1-\chi(s/a)),\quad s\in \R.$$
It is easy to see that $\Sigma=\varphi_\Sigma^{-1}(0)$, and $\nabla \varphi_\Sigma(x)=\nu_\Sigma(x)$, for each $x\in \Sigma$. Moreover, $\kappa=0$ is an eigenvalue of $\nabla^2\varphi_\Sigma(x)$, the remaining eigenvalues of $\nabla^2\varphi_\Sigma(x)$ are the principal curvatures $\kappa_j$ of $\Sigma$ at $x\in\Sigma$.

If we consider the subset $\cMH^2(\Omega,r)$ of $\cMH^2(\Omega)$ which consists of all closed hypersurfaces $\Gamma\in \cMH^2(\Omega)$ such that $\Gamma\subset \Omega$ satisfies the ball condition with fixed radius $r>0$
(see the definition below) then the map $\Phi:\cMH^2(\Omega,r)\to C^2(\bar{\Omega})$ defined by $\Phi(\Gamma)=\varphi_\Gamma$ is an isomorphism of
the metric space $\cMH^2(\Omega,r)$ onto $\Phi(\cMH^2(\Omega,r))\subset C^2(\bar{\Omega})$.

Let $s-(n-1)/p>2$; for $\Gamma\in\cMH^2(\Omega,r)$, we define $\Gamma\in W^s_p(\Omega,r)$ if $\varphi_\Gamma\in W^s_p(\Omega)$. In this case the local charts for $\Gamma$ can be chosen of class $W^s_p$ as well. A subset $A\subset W^s_p(\Omega,r)$ is  (relatively) compact, if and only if  $\Phi(A)\subset W^s_p(\Omega)$ is (relatively) compact.

As an ambient space for the
state manifold $\cSM$ of problem \eqref{i2pp}  we consider
the product space $C(\bar{\Omega})^{n+1}\times \cMH^2(\Omega)$, due to continuity of velocity, temperature
and curvature.

We define the state manifold $\cSM$ as follows.
\begin{equation*}
\begin{aligned}
\label{phasemanif0}
\cSM:=&
\Big\{(u,\theta,\Gamma)\in C(\bar{\Omega})^{n+1}\times \cMH^2: (u,\theta)\in W^{2-2/p}_p(\Omega\setminus\Gamma)^{n+1},\, \Gamma\in W^{4-3/p}_p, \\
 & {\rm div}\, u=0\; \mbox{ in }\; \Omega\setminus\Gamma,\;\ \theta>0 \mbox{ in } \bar{\Omega},
 \;\ u=\partial_\nu\theta =0\; \mbox{ on } \partial\Omega,\\
& [\![\psi(\theta)]\!]+\sigma H_\Gamma=P_\Gamma[\![\mu(\theta)(\nabla u+[\nabla u]^{\sf T}) \nu_\Gamma]\!] =0\mbox{ on } \Gamma,
\;\ l(\theta)\neq0\mbox{ on } \Gamma,\;\\
&[\![d(\theta)\partial_\nu \theta]\!] \in W^{2-6/p}_p(\Gamma)\Big\},
\end{aligned}
\end{equation*}
Charts for these manifolds are obtained by the charts induced by $\cMH^2(\Omega)$,
followed by a Hanzawa transformation as defined in Section 4.

Applying Theorem \ref{wellposed} and re-parameterizing the interface repeatedly,
we see that (\ref{i2pp}) yields a local semiflow on $\cSM$.

\bigskip
%%%%%%%%%%%%%%%%
\begin{theorem}
\label{semiflow} Let $p>n+2$, $\sigma>0$, and suppose $\psi_i\in C^3(0,\infty)$, $\mu_i,d_i\in C^2(0,\infty)$ such that
$$\kappa_i(s)=-s\psi_i^{\prime\prime}(s)>0,\quad \mu_i(s)>0,\quad  d_i(s)>0,\quad s\in(0,\infty),\; i=1,2.$$
Then problem (\ref{i2pp}) generates a local semiflow
on the state manifold $\cSM$. Each solution $(u,\theta,\Gamma)$ exists on a maximal time
interval $[0,t_*)$, where $t_*=t_*(u_0,\theta_0,\Gamma_0)$.
\end{theorem}
%%%%%%%%%%%%%%%%%%%%%
\bigskip

\noindent
Note that the pressure does not occur explicitly as a variable in the
local semiflow, as the latter is only formulated in terms of the temperature $\theta$,
the velocity field $u$, and the free boundary $\Gamma$.
The pressure $\pi$ is determined for each $t$ from
$(u,\theta,\Gamma)$ by means of the weak transmission problem
\begin{align*}
(\nabla\pi|\nabla\phi)_{L_2(\Omega)} &= 2({\rm div}(\mu(\theta)D)-u\cdot\nabla u|\nabla\phi)_{L_2(\Omega)},\quad \phi\in H^1_{p^\prime}(\Omega),\\
[\![\pi]\!]&= 2[\![\mu(\theta)D]\!]\quad \mbox{ on } \Gamma.
\end{align*}
Concerning such transmission problems we refer to Section 8 in \cite{KPW10}.

\bigskip

\noindent
Let $(u,\theta,\Gamma)$ be a solution in the state manifold $\cSM$ with maximal interval of existence $[0,t_*)$. By the
{\em uniform ball condition} we mean the existence of a radius $r_0>0$ such that for each $t\in[0,t_*)$,
at each point $x\in\Gamma(t)$ there exists centers $x_i\in \Omega_i(t)$ such that
$B_{r_0}(x_i)\subset \Omega_i$ and $\Gamma(t)\cap \bar{B}_{r_0}(x_i)=\{x\}$, $i=1,2$. Note that this condition
bounds the curvature of $\Gamma(t)$, prevents parts of it to shrink to points, to touch the outer
boundary $\partial \Omega$, and to undergo topological changes.

With this property, combining the local semiflow for (\ref{i2pp}) with
the Ljapunov functional and compactness we obtain the following result.

\bigskip

\begin{theorem} \label{Qual} Let $p>n+2$, $\sigma>0$, and suppose $\psi_i\in C^3(0,\infty)$, $\mu_i,d_i\in C^2(0,\infty)$ such that
$$\kappa_i(s)=-s\psi_i^{\prime\prime}(s)>0,\quad \mu_i(s)>0,\quad  d_i(s)>0,\quad s\in(0,\infty),\; i=1,2.$$
Suppose that $(u,\theta,\Gamma)$ is a solution of
(\ref{i2pp}) in the state manifold $\cSM$ on its maximal time interval $[0,t_*)$.
Assume there is constant $M>0$  such that the  following conditions hold on $[0,t_*)$:
\begin{itemize}

\item[(i)] $|u(t)|_{{W^{2-2/p}_p}^n},|\theta(t)|_{W^{2-2/p}_p},|\Gamma(t)|_{W^{4-3/p}_p},
|[\![d(\theta(t))\partial_\nu u(t)]\!]|_{W^{2-6/p}_p}\leq M\!<\!\infty$;

\item[(ii)] $|l(\theta(t))|, \theta(t)\geq 1/M$;

\item[(iii)]  $\Gamma(t)$ satisfies the uniform ball condition.
\end{itemize}

\noindent
Then $t_*=\infty$, i.e.\ the solution exists globally, and  $\omega_+(u,\theta,\Gamma)\subset\cE$ is non-empty.
\end{theorem}

\begin{proof}
Assume that  (i), (ii), and (iii)  are valid. Then $\Gamma([0,t_*))\subset W^{4-3/p}_p(\Omega,r)$ is bounded, hence relatively compact in
$W^{4-3/p-\ve}_p(\Omega,r)$. Thus we may cover this set by finitely many balls with centers $\Sigma_k$ real analytic in such a way that
${\rm dist}_{W^{4-3/p-\ve}_p}(\Gamma(t),\Sigma_j)\leq \delta$ for some $j=j(t)$, $t\in[0,t_*)$. Let $J_k=\{t\in[0,t_*):\, j(t)=k\}$; using for each $k$ a Hanzawa-transformation $\Xi_k$, we see that the pull backs $\{(u(t,\cdot),\theta(t,\cdot))\circ\Xi_k:\, t\in J_k\}$ are bounded in $W^{2-2/p}_p(\Omega\setminus \Sigma_k)^{n+1}$, hence relatively compact in $W^{2-2/p-\ve}_p(\Omega\setminus\Sigma_k)^{n+1}$. Employing now Corollary \ref{wellposed3}  we obtain solutions
 $(u^1,\theta^1,\Gamma^1)$ with initial configurations $(u(t),\theta(t),\Gamma(t))$ in the state manifold on a common time interval, say $(0,\tau]$, and by uniqueness we have
$$(u^1(\tau),\theta^1(\tau),\Gamma^1(a))=(u(t+\tau),\theta(t+\tau),\Gamma(t+\tau)).$$
 Continuous dependence implies that the  orbit of the solution $(u(\cdot),\theta(\cdot),\Gamma(\cdot))$ is  relative compact in $\cSM$, in particular $t_*=\infty$ and  $(u,\theta,\Gamma)(\R_+)\subset\cSM$ is relatively compact.
The negative total entropy is a strict Ljapunov functional, hence the limit set $\omega_+(u,\theta,\Gamma)\subset \cSM$
of a solution is contained in the set $\cE$ of equilibria.  By compactness $\omega_+(u,\theta,\Gamma)\subset \cSM$ is non-empty, hence the solution comes close to $\cE$, and stays there.
\end{proof}

\bigskip

\noindent
{\bf Remark.}\, We can prove that each equilibrium is normally hyperbolic. Therefore, taking into account the results in \cite{PSZ10}, one may expect that each orbit satisfying the assumptions of
Theorem~\ref{Qual} actually converges to an equilibrium. This topic will be considered elsewhere.

\bigskip

{
}

\end{document}